\documentclass[10pt]{amsart}
\usepackage{amsfonts}
\usepackage{amssymb,latexsym}
\usepackage{amsmath}
\usepackage{amsthm}
\usepackage{amssymb}
\usepackage{enumerate}
\newtheorem{theorem}{Theorem}[section]
\newtheorem{lemma}[theorem]{Lemma}
\newtheorem{proposition}[theorem]{Proposition}
\newtheorem{corollary}[theorem]{Corollary}
\theoremstyle{definition}
\newtheorem{definition}[theorem]{Definition}

\theoremstyle{remark}
\newtheorem{remark}[theorem]{Remark}
\numberwithin{equation}{section}

\begin{document}
\setlength{\baselineskip}{1.2\baselineskip}
\title  [A CONSTANT RANK THEOREM FOR PARTIAL CONVEX SOLUTIONS]
{A CONSTANT RANK THEOREM FOR PARTIAL CONVEX SOLUTIONS OF PARTIAL
DIFFERENTIAL EQUATIONS}
\author{Chuanqiang Chen}
\address{Department of Mathematics\\
         University of Science and Technology of China\\
         Hefei, 230026, Anhui Province, CHINA}
\email{cqchen@mail.ustc.edu.cn}
\thanks{Research of the author was supported by Grant 10871187 from the National Natural Science
Foundation of China.}
\maketitle

\begin{abstract}
Thanks to the test function of Bian-Guan[2], we successfully obtain
a constant rank theorem for partial convex solutions of a class
partial differential equations. This is the microscopic version of
the macroscopic partial convexity principle in [1], and also is a
generalization of the result in [2].
\end{abstract}

\section{Introduction}
The convex solution of partial differential equation is an
interesting issue for a long time. And so far as we know, there
are two important methods for this problem, which are macroscopic
and microscopic methods. Whereas there are many solutions which
are not convex. For example, the admissible solutions of the
Hessian equations were studied in [7,10], the power concave
solutions in [13,18], and the $k$-convex solutions in [11]. In
this paper we will consider the partial convex solutions (see [1]
or Definition 1.1 as below) of the elliptic and parabolic
equations.

The study of macroscopic convexity is using a weak maximum
principle, while the study of microscopic convexity is using a
strong maximum principle. For the macroscopic convexity argument,
Korevaar made breakthroughs in [14,15], he introduced concavity
maximum principles for a class of quasilinear elliptic equations.
And later it was improved by Kennington [13] and by Kawhol [12].
The theory further developed to its great generality by
Alvarez-Lasry-Lions [1]. The key of the study of microscopic
convexity is a method called \emph{constant rank theorem} which
was discovered in 2 dimension by Caffarelli-Friedman [5] (a
similar result was also discovered by Singer-Wong-Yau-Yau [19] at
the same time). Later the result in [5] was generalized to
$\mathbb{R}^n$ by Korevaar-Lewis [17]. Recently the \emph{constant
rank theorem} was generalized to fully nonlinear equations in [6]
and [2], where the result in [2] is the microscopic version of the
macroscopic convexity principle in [1].

\emph{Constant rank theorem} is a very useful tool to produce convex
solutions in geometric analysis. By the corresponding homotopic
deformation, the existence of convex solution comes from the
\emph{constant rank theorem}. For the geometric application of the
\emph{constant rank theorem}, the Christoffel-Minkowski problem and
the related prescribing Weingarten curvature problems were studied
in [8,9,10]. The preservation of convexity for the general geometric
flows of hypersurfaces was given in [2]. Soon after the
\emph{constant rank theorem} for the level set was established in
[3], where [3] is a microscopic version of [4] (also it was studied
in [16]). And the existence of the $k$-convex hypersurface with
prescribed mean curvature was given in [11] recently.

In this paper we consider the partial convexity of solutions of the
following elliptic equation, and give a constant rank theorem for
partial convex solutions
\begin{equation}
F(D^2 u,Du,u,x) = 0, \quad x \in \Omega  \subset \mathbb{R}^N,
\end{equation}
where $ F \in C^{2,1} (\mathcal{S}^N \times \mathbb{R}^N \times
\mathbb{R} \times \Omega )$ and $F$ is elliptic in the following
sense
\begin{equation}
(\frac{{\partial F}} {{\partial u_{ab} }}(D^2 u,Du,u,x))_{N \times
N}  > 0,  \quad \text{for all } x \in \Omega.
\end{equation}

First, we give the definition of the partial convexity of a function
$u$, which could be found in [1].
\begin{definition}
Suppose $u \in C^2 (\Omega ) \cap C(\overline\Omega  )$, where $
\Omega$ is a domain in $\mathbb{R}^N  = \mathbb{R}^{N'} \times
\mathbb{R}^{N'' } $, and $N'$ and $N''$are two integers with
$N=N'+N''$. Then $u$ is partial convex (with respect to the first
variable) is that $ x' \to u(x',x'') $ is convex for every
$x=(x',x'') \in \overline \Omega$. In particular, if $N''=0$, i.e.
$u$ is convex in $\Omega \subset \mathbb{R}^N=\mathbb{R}^{N'}$, $u$
is said degenerate partial convex.
\end{definition}
For simplicity, we introduce additional notations . As in [1], we
denote $\mathcal {S}^n$ to be the set of all real symmetric $n
\times n$ matrices. And we shall write $p \in \mathbb{R}^N$ in the
form $(p',p'')$ with $p' \in \mathbb{R}^{N'}$, $p'' \in
\mathbb{R}^{N''}$ and split a matrix $A \in \mathcal{S}^N$ into $
\left( {\begin{matrix}
   a & b  \\
   {b^T } & c  \\
 \end{matrix} } \right)
$ with $a \in \mathcal{S}^{N'}$, $b \in \mathbb{R}^{N' \times N''}$
and $c \in \mathcal{S}^{N''}$; we also let
\begin{equation} F(A,p,u,x)=F(\left(
{\begin{matrix}
   a & b  \\
   {b^T } & c  \\
 \end{matrix} } \right),p',p'',u,x',x'').
\end{equation}

One of our main results is the following theorem
\begin{theorem}\label{TH}(CONSTANT RANK THEOREM)
Suppose $\Omega$ is a domain in $\mathbb{R}^N  = \mathbb{R}^{N'}
\times \mathbb{R}^{N'' } $ and $ F(A,p,u,x) \in C^{2,1}
(\mathcal{S}^N \times \mathbb{R}^N \times \mathbb{R} \times \Omega
)$. If $F$ satisfies (1.2) and the following condition
\begin{equation}
F(\left( {\begin{matrix}
   a^{-1} & a^{-1}b  \\
   {({a^{-1}b})^T } & c+ b^Ta^{-1}b \\
 \end{matrix} } \right),p',p'',u,x',x'')  \text{ is locally convex in } (a,b,c,p'',u,x').
\end{equation}
If $u \in C^{2,1} (\Omega)$ is a partial convex solution of (1.1),
then $(u_{ij})_{N' \times N'}$ has constant rank in $\Omega$.
\end{theorem}

\begin{remark}
if $N''=0$, i.e. for the degenerate partial convexity, structure
condition (1.4) is \emph{inverse-convex}  condition, the result of
Bian-Guan [2]. And for general partial convexity, structure
condition (1.4) is strictly stronger than the \emph{inverse-convex}
condition.
\end{remark}

An immediate consequence of Theorem 1.2 is for partial convex
solutions of the following quasilinear second elliptic equation
\begin{equation}
 \sum\limits_{a,b = 1}^N
{a^{ab} (x'',u_1 (x), \cdots ,u_{N'} (x))u_{ab} (x)}  =
f(x,u(x),Du(x)) > 0,
\end{equation}
where $x \in \Omega  \subset \mathbb{R}^N $ and
\begin{equation}
(a^{ab} (x'',u_1 (x), \cdots ,u_{N'} (x)))_{N \times N}  > 0,  \quad
\text{for all } x \in \Omega.
\end{equation}

\begin{corollary}
Suppose $\Omega$ is a domain in $\mathbb{R}^N  = \mathbb{R}^{N'}
\times \mathbb{R}^{N'' } $, and $u \in C^{2,1} (\Omega) $ is the
partial convex solution of (1.5). If
\begin{center}
$f(x',x'',u,p',p'')$ is locally concave in $(p'',u,x')$,
\end{center}
then $(u_{ij})_{N' \times N'}$ has constant rank in $\Omega$.
\end{corollary}
Set
\begin{equation}
F(D^2 u,Du,u,x)=\sum\limits_{a,b = 1}^N {a^{ab} (x'',u_1 (x), \cdots
,u_{N'} (x))u_{ab} (x)}-f(x,u(x),Du(x)),
\end{equation}
we can verify that $F$ satisfies the structure condition (1.4) (see
the equivalent condition (3.13) in the third section).

A corresponding result holds for the parabolic equation.
\begin{theorem}
Suppose $\Omega$ is a domain in $\mathbb{R}^N  = \mathbb{R}^{N'}
\times \mathbb{R}^{N'' } $, and $ F(A,p,u,x,t) \in C^{2,1}
(\mathcal{S}^N \times \mathbb{R}^N \times \mathbb{R} \times \Omega
\times (0,T])$. If $F$ satisfies (1.2) for each $t$ and the
following condition
\begin{equation}
F(\left( {\begin{matrix}
   a^{-1} & a^{-1}b  \\
   {({a^{-1}b})^T } & c+ b^Ta^{-1}b \\
 \end{matrix} } \right),p',p'',u,x',x'',t)  \text{ is locally convex in } (a,b,c,p'',u,x').
\end{equation}
If $u \in C^{2,1} (\Omega \times (0,T]) $ is a partial convex
solution of the equation
\begin{equation}
\frac{{\partial u}} {{\partial t}} = F(D^2 u,Du,u,x,t), \quad (x,t)
\in \Omega \times (0,T],
\end{equation}
then $(u_{ij}(x,t))_{N' \times N'}$ has constant rank in $\Omega$
for each $T \geqslant t >0$. Moreover, let $l(t)$ be the minimal
rank of $(u_{ij}(x,t))_{N' \times N'}$ in $\Omega$, then $l(s)
\leqslant l(t)$ for all $s \leqslant t \leqslant T$.
\end{theorem}
The rest of the paper is organized as follows. In section 2, we
work on the Laplace equation, a special case of Corollary 1.4. In
section 3, using the key auxiliary function $q(x)$ in [2], we do
some preliminarily calculations on the constant rank theorem. In
section 4, we prove the Theorem 1.2 using a strong maximum
principle. In section 5, we prove Theorem 1.5. And the last
section is devoted to a discussion of the structure condition.

\textbf{Acknowledgement}. The author would like to express sincere
gratitude to Prof. Xi-Nan Ma for his encouragement and many
suggestions in this subject.

\section{an example}

In this section, we give a constant rank theorem for partial
convex solutions of Laplace equation, a special case of Corollary
1.4.

We rewrite the result as follows.
\begin{theorem}
Suppose $\Omega$ is a domain in $\mathbb{R}^N  = \mathbb{R}^{N'}
\times \mathbb{R}^{N'' } $, and $u \in C^{2,1} (\Omega)$ is the
partial convex solution of the following equation
\begin{equation}
 \Delta u(x)=\sum\limits_{a = 1}^N{ u_{aa} (x)}  = f(x,u(x),Du(x)) >
 0, \quad x \in \Omega.
\end{equation}
Assume
\begin{align}
f(x',x'',u,p',p'') \text{ is locally concave in }(p'',u,x'),
\end{align}
then $(u_{ij})_{N' \times N'}$ has constant rank in $\Omega$.
\end{theorem}
Before the proof of Theorem 2.1, we do some preliminaries. As in
[9], we recall the definition of $k$-symmetric functions: For $1
\leqslant k \leqslant N'$, and $ \lambda  = (\lambda _1 ,\lambda _2
, \cdots ,\lambda _{N'} ) \in \mathbb{R}^{N'} $,
$$
\sigma _k (\lambda ) = \sum\limits_{i_1  < i_2  <  \cdots  < i_k }
{\lambda _{i_1 } \lambda _{i_2 }  \cdots \lambda _{i_k } },
$$
we denote by $\sigma _k (\lambda \left| i \right.)$ the symmetric
function with $\lambda_i = 0$ and $\sigma _k (\lambda \left| ij
\right.)$ the symmetric function with $\lambda_i =\lambda_j = 0$.

The definition can be extended to symmetric matrices by letting
$\sigma_k(W) = \sigma_k(\lambda(W))$, where $ \lambda(W)= (\lambda
_1(W),\lambda _2 (W), \cdots ,\lambda _{N'}(W))$ are the eigenvalues
of the symmetric matrix $W$. We also set $\sigma_0 = 1$ and
$\sigma_k = 0$ for $k > N'$.

We need the following standard formulas, which could be found in
[9], [2] or [3].
\begin{lemma}
Suppose $W=(W_{ij})$ is diagonal, and $m$ is positive integer, then
$$
\frac{{\partial \sigma _m (W)}} {{\partial W_{ij} }} = \left\{
\begin{matrix}
\sigma _{m - 1} (W\left| i \right.), \quad if \quad i = j, \hfill\cr
  0, \qquad if \quad i \ne j. \hfill \cr
 \end{matrix}  \right.
$$

$$
\frac{{\partial ^2 \sigma _m (W)}} {{\partial W_{ij} \partial W_{kl}
}} = \left\{ \begin{matrix}
  \sigma _{m - 2} (W\left| {ik} \right.), \quad if \quad i = j,k = l,i \ne k, \hfill \cr
   - \sigma _{m - 2} (W\left| {ik} \right.), \quad if \quad i = l,j = k,i \ne j, \hfill \cr
  0, \qquad \qquad otherwise. \hfill \cr
 \end{matrix}  \right.
$$
\end{lemma}

\textbf{Proof of Theorem 2.1}. With the assumptions of $u$ in
Theorem 2.1, $u$ is automatically in $C^{3,1}$. We denote
$W=(u_{ij})_{N' \times N'}$. For each $z_0 \in \Omega$ where $W$ is
of minimal rank $l$. We pick a small open neighborhood $\mathcal
{O}$ of $z_0$, we will prove it always be rank of $l$ in $\mathcal
{O}$.
   We shall use the strong minimum principle to prove the theorem.
Let
\begin{align}
 \phi (x)= \sigma _{l + 1} (W),
\end{align}
then $ \phi(z_0)=0 $. We shall show $\phi(x)\equiv 0$ in $\mathcal
{O}$. If true, it implies  the set $ \left\{ {x \in \Omega|\phi (x)
= 0} \right\} $ is an open set. But it is also closed, then we get $
\phi (x) \equiv 0$ in $\Omega$ since $\Omega$ connected, i.e.
$(u_{ij})_{N' \times N'}$ is of constant rank $l$ in $\Omega$.

 Following Caffarelli and Friedman [5], for two functions $h(y)$ and $k(y)$
defined in an open set $ \mathcal {O} \subset \Omega$, we say that
$h(y)\lesssim k(y)$ provided there exist positive constants $c_1$
and $c_2$ such that
\begin{equation}
(h-k)(y)\leq (c_1 |\nabla \phi| + c_2 \phi)(y).
\end{equation}
We also write  $h(y)\sim k(y)$  if $h(y)\lesssim k(y)$ and
$k(y)\lesssim h(y)$ . Next, we write $h\lesssim k$ if the above
inequality holds in the neighborhood $\mathcal {O}$, with the
constants $c_1$ and $c_2$ independent of $y$ in this neighborhood.
Finally $h\sim k$ if $h\lesssim k$ and $k\lesssim h$.

We shall show that
\begin{equation}
\Delta \phi(x) =\sum\limits_{a = 1}^N { \phi _{aa} (x)}\mathop
\lesssim 0.
\end{equation}
Since $ {\phi (x) \geqslant 0} $ in $\Omega$ and $\phi(z_0)=0$, it
then follows from the Strong Minimum Principle that $\phi(x) \equiv
0$ in $\mathcal {O}$.

For any fixed point $x \in \mathcal {O}$, we rotate coordinate $e_1,
\cdots, e_{N'}$ such that the matrix ${u_{ij}},i,j=1, \cdots, N'$ is
diagonal and without loss of generality we assume $ u_{11} \leqslant
u_{22} \leqslant \cdots \leqslant u_{N'N'} $. Then there is a
positive constant $C > 0$ depending only on $\left\| u
\right\|_{C^{3,1} }$ and $\mathcal {O}$, such that $ u_{N'N'}
\geqslant \cdots \geqslant u_{N' - l+1N' - l+1} \geqslant C > 0 $
for all $x \in \mathcal {O}$. For convenience we denote $ G = \{ N'
- l+1, \cdots ,N'\} $ and $ B = \{ 1,2, \cdots ,N'- l\} $ which
means good terms and bad ones in indices respectively. Without
confusion we will also simply denote $ B = \{ u_{11} , \cdots
,u_{N'-lN'-l} \} $ and $ G = \{ u_{N' - l+1N' - l+1} , \cdots
,u_{N'N'} \} $. In the following, all the calculation at the point
$x$ are using the relation $\lesssim$ with the understanding that
the constants in (2.4) are under control.

Following a direct computation as in [9] and $W$ is diagonal, we can
get
\begin{align}
&0 \sim \phi \sim \sigma _l (G)\sum\limits_{i \in B} {u_{ii} } ,
\text { and } u_{ii} \sim 0 \text{ for each } i \in B;\\
&0 \sim \phi _a   \sim \sigma _l (G)\sum\limits_{i \in B} {u_{iia
}};
\end{align}
then by (2.6),(2.7) and Lemma 2.2, we obtain
\begin{align}
&\Delta \phi  = \sum\limits_{a = 1}^N {\frac{{\partial ^2 \phi }}
{{\partial x_a \partial x_a }}}  = \sum\limits_{a = 1}^N
{[\sum\limits_{i,j = 1}^{N'} {\frac{{\partial \sigma _{l + 1} (W)}}
{{\partial u_{ij} }}u_{ijaa}  + \sum\limits_{i,j,k,l = 1}^{N'}
{\frac{{\partial ^2 \sigma _{l + 1} (W)}} {{\partial u_{ij} \partial
u_{kl} }}u_{ija} u_{kla} } ]} }  \notag \\
&= \sum\limits_{a = 1}^N {[\sum\limits_{i = 1}^{N'} {\frac{{\partial
\sigma _{l + 1} (W)}} {{\partial u_{ii} }}u_{iiaa}  +
\sum\limits_{i,j = 1}^{N'} {\frac{{\partial ^2 \sigma _{l + 1} (W)}}
{{\partial u_{ii} \partial u_{jj} }}u_{iia} u_{jja}  +
\sum\limits_{i,j = 1}^{N'} {\frac{{\partial ^2 \sigma _{l + 1} (W)}}
{{\partial u_{ij} \partial u_{ji} }}u_{ija} u_{jia} } ]} } } \notag \\
& \sim \sum\limits_{a = 1}^N {[\sigma _l (G)\sum\limits_{i \in B}
{u_{iiaa} }  - 2\sigma _l (G)\sum\limits_{i \in B,j \in G} {\frac{1}
{{u_{jj} }}u_{ija} u_{jia} } ]} \\
 &\sim \sigma _l (G)\sum\limits_{i \in B} {(\Delta u)_{ii} }
 - 2\sigma _l (G)\sum\limits_{a = 1}^N {\sum\limits_{i \in B,j \in G} {\frac{1}
{{u_{jj} }}u_{ija} ^2 } }. \notag
\end{align}
For each $i \in B$, we differentiate (2.1) twice in $x_i$, then we
obtain
\begin{align*}
(\Delta u)_{ii}=&[f_{x_i}+f_u u_i+\sum\limits_{a = 1}^N
f_{p_a} u_{a i} ]_i \\
=&f_{x_i x_i }  + 2f_{u,x_i } u_i  + f_{u,u} u_i ^2 \\
&+ 2\sum\limits_{a = 1}^N {f_{x_i ,p_a } u_{a i} } + 2\sum\limits_{a
= 1}^N {f_{u,p_a  } u_i u_{a i} } + \sum\limits_{a ,b = 1}^N
{f_{p_a,p_b } u_{a i} u_{b i} },
\end{align*}
since $W=(u_{ij})_{N' \times N'}$ is diagonal and (2.6), we get from
the above equation
\begin{align}
(\Delta u)_{ii} \sim& f_{x_i x_i }  + 2f_{u,x_i } u_i  + f_{u,u} u_i ^2 \notag \\
&+ 2\sum\limits_{\alpha = N' + 1}^N {f_{x_i ,p_\alpha  } u_{\alpha
i} } + 2\sum\limits_{\alpha = N' + 1}^N {f_{u,p_\alpha  } u_i
u_{\alpha i} } + \sum\limits_{\alpha ,\beta  = N' + 1}^N
{f_{p_\alpha ,p_\beta } u_{\alpha i} u_{\beta i} },
\end{align}
so we obtain from (2.8) and (2.9)
\begin{align}
\frac{\Delta \phi}{\sigma _l (G)} \sim &\sum\limits_{i =1}^{N'-l}
{(\Delta u)_{ii} }  - 2\sum\limits_{a = 1}^N {\sum\limits_{j
=N'-l+1}^{N'}{\frac{1} {{u_{jj} }}\sum\limits_{i = 1}^{N' - l} {u_{ija} ^2 } } } \notag \\
\sim &- 2\sum\limits_{a = 1}^N {\sum\limits_{j =N'-l+1}^{N'}
{\frac{1} {{u_{jj} }}\sum\limits_{i = 1}^{N' - l} {u_{ija} ^2 } } }
+\sum\limits_{i =1}^{N'-l} {[ f_{x_i x_i }  + 2f_{u,x_i } u_i  +
f_{u,u} u_i ^2 } \\
&+ 2\sum\limits_{\alpha = N' + 1}^N {f_{x_i ,p_\alpha } u_{\alpha i}
} + 2\sum\limits_{\alpha = N' + 1}^N {f_{u,p_\alpha } u_i u_{\alpha
i} } + \sum\limits_{\alpha ,\beta  = N' + 1}^N {f_{p_\alpha ,p_\beta
} u_{\alpha i} u_{\beta i} } ].  \notag
\end{align}
 By the condition (2.2), we obtain (2.5).  The proof of Theorem
2.1 is completed.
\begin{remark}
In (2.8), we have used Lemma 2.5 in [2], otherwise the first
"$\sim$" will be "$\lesssim$".
\end{remark}
\begin{remark}
By a similar proof as above, we can get the general case of
Corollary 1.4.
\end{remark}

\section{primarily calculations on the constant rank theorem}

\subsection{calculations on the test function}

With the assumptions in Theorem 1.2 and Theorem 1.5, $u$ is in
$C^{3,1}$. Let $ W= (u_{ij} )_{N' \times N'} $ and $ l = \mathop
{\min }\limits_{x \in \Omega } rank(W(x)) $. We may assume $l
\leqslant N' - 1 $, otherwise there is nothing to prove. Suppose
$z_0 \in \Omega$ is a point where $W$ is of minimal rank $l$.

Throughout this paper we assume that $ 1 \leqslant i,j,k,l,m,n
\leqslant N' $, $ N' \leqslant \alpha ,\beta ,\gamma ,\eta ,\xi
,\zeta  \leqslant N'' $, $ 1 \leqslant a,b,c,d \leqslant N $ and $
\sigma_j(W) = 0$ if $j < 0$ or $j > N'$. As in Bian-Guan [2],we
define for $W=(u_{ij}(x)) \in \mathcal{S}^{N'}$,
\begin{equation}
q(W) = \left\{ \begin{matrix}
  \frac{{\sigma _{l + 2} (W)}}{{\sigma _{l + 1} (W)}}, \quad if \quad\sigma _{l + 1} (W) > 0, \hfill \cr
  0, \qquad if \quad \sigma _{l + 1} (W) = 0. \hfill \cr
 \end{matrix}  \right.
\end{equation}
and we consider the following test function
\begin{align}
 \phi = \sigma _{l + 1} (W) + q(W).
\end{align}

For each $z_0 \in \Omega$ where $W$ is of minimal rank $l$. We pick
an open neighborhood $\mathcal {O}$ of $z_0$, and for any fixed
point $x \in \mathcal {O}$, we rotate coordinate $e_1, \cdots,
e_{N'}$ such that the matrix ${u_{ij}},i,j=1, \cdots, N'$ is
diagonal and without loss of generality we assume $ u_{11} \leqslant
u_{22} \leqslant \cdots \leqslant u_{N'N'} $. Then there is a
positive constant $C > 0$ depending only on $\left\| u
\right\|_{C^{3,1} }$ and $\mathcal {O}$, such that $ u_{N'N'}
\geqslant \cdots \geqslant u_{N' - l+1N' - l+1} \geqslant C > 0 $
for all $x \in \mathcal {O}$. For convenience we denote $ G = \{ N'
- l+1, \cdots ,N'\} $ and $ B = \{ 1,2, \cdots ,N'- l\} $ which
means good terms and bad ones in indices respectively. Without
confusion we will also simply denote $ B = \{ u_{11} , \cdots
,u_{N'-lN'-l} \} $ and $ G = \{ u_{N' - l+1N' - l+1} , \cdots
,u_{N'N'} \} $. Note that for any $\delta > 0$, we may choose
$\mathcal {O}$ small enough such that $u_{jj} < \delta$ for all $j
\in B$ and $x \in \mathcal {O}$.

We will use notation $h = O(f)$ if $ \left| {h(x)} \right| \leqslant
Cf(x)$ for $x \in \mathcal {O}$ with positive constant $C$ under
control. It is clear that $u_{ii} = O(\phi)$ for all $i \in B$.

To get around $\sigma_{l+1}(W) = 0$, for $\varepsilon> 0 $
sufficient small, we consider
\begin{align}
q(W_\varepsilon )=\frac{{\sigma _{l + 2} (W_\varepsilon)}}{{\sigma
_{l + 1} (W_\varepsilon)}}, \quad \phi _\varepsilon  = \sigma _{l +
1} (W_\varepsilon ) + q(W_\varepsilon ),
\end{align}
where $W_\varepsilon  = W + \varepsilon I$. We will also denote $
B_\varepsilon   = \{ u_{11}  + \varepsilon , \cdots ,u_{N'- lN'-l} +
\varepsilon \} $, $ G_\varepsilon   = \{ u_{N' - l+1N' - l+1} +
\varepsilon , \cdots ,u_{N'N'}  + \varepsilon \} $. (see
Bian-Guan[2]).

Set $u_\varepsilon (x)=u(x)+\frac{\varepsilon}{2} \left| {x'}
\right|^2 $, then $ W_\varepsilon= ((u_\varepsilon)_{ij} )_{N'
\times N'} $. To simplify the notations, we will write $u$ for
$u_\varepsilon$, $q$ for $q_\varepsilon$, $W$ for $W_\varepsilon$,
$G$ for $G_\varepsilon$, and $B$ for $B_\varepsilon$ with the
understanding that all the estimates will be independent of
$\varepsilon$. In this setting, if we pick $\mathcal {O}$ small
enough, there is $C > 0$ independent of $\varepsilon$ such that
\begin{equation}
\phi  \geqslant C\varepsilon ,\sigma _1 (B) \geqslant C\varepsilon ,
 \text{ for all } x\in \mathcal {O}.
\end{equation}

First, we consider the regularity of $q(W(x))$.
\begin{proposition} ([2])
let $u \in C^{3,1} (\Omega) $ be a partial convex function with the
first variable and $ W(x)= (u_{ij}(x) )_{N' \times N'} $. Let $ l =
\mathop {\min }\limits_{x \in \Omega } rank(W(x)) $, then the
function $q(x) = q(W(x))$ defined in (3.1) is in $C^{1,1} (\Omega)
$.
\end{proposition}

In the following, we denote
\begin{align*}
& F^{ab}  = \frac{{\partial F}} {{\partial u_{ab} }} ,  F^{p_a } =
\frac{{\partial F}} {{\partial u_a }} ,  F^u  = \frac{{\partial F}}
{{\partial u}} ,\\
&F^{ab,cd}  = \frac{{\partial ^2 F}} {{\partial u_{ab} \partial
u_{cd} }},  F^{ab,p_c }  = \frac{{\partial ^2 F}} {{\partial
u_{ab}\partial u_c }} ,  F^{ab,u}  = \frac{{\partial ^2 F}}
{{\partial u_{ab}\partial u}},\\
& F^{p_a p_b }  = \frac{{\partial ^2 F}} {{\partial u_a \partial
u_b}} ,  F^{p_a, u}  = \frac{{\partial ^2 F}} {{\partial u_a
\partial u}} ,  F^{u,u}  = \frac{{\partial ^2 F}} {{\partial
u\partial u}},
\end{align*}
where $ 1 \leqslant a,b,c \leqslant N $.

\begin{theorem}
Suppose $\Omega$ is a domain in $\mathbb{R}^N  = \mathbb{R}^{N'}
\times \mathbb{R}^{N'' } $ and $u \in C^{3,1} (\Omega) $ is a
partial convex solution of (1.1). Let $l$ be the minimal rank of $
W= (u_{ij} )_{N' \times N'} $ in $\Omega$. Suppose $l$ is attained
at $z_0 \in \Omega$, and $\mathcal {O}$ is a small neighborhood of
$z_0$ as above. For any fixed point $x \in \mathcal {O}$ we choose
the coordinate such that $W(x)$ is diagonal. Then at $x$ we have
\begin{eqnarray}
&\sum\limits_{a,b = 1}^N {F^{ab} \phi _{ab} }&  = \sum\limits_{i \in
B} {[\sigma _l (G) + \frac{{\sigma _1 ^2 (B\left| i \right.) -
\sigma _2 (B\left| i \right.)}} {{\sigma _1 ^2
(B)}}]\sum\limits_{a,b= 1}^N {F^{ab} u_{iiab} } } \notag \\
&&- 2\sum\limits_{i \in B,j \in G} {[\sigma _l (G) + \frac{{\sigma
_1 ^2 (B\left| i \right.) - \sigma _2 (B\left| i \right.)}} {{\sigma
_1 ^2 (B)}}]\frac{1} {{u_{jj} }}\sum\limits_{a,b = 1}^N {F^{ab}
u_{ija}u_{ijb} } } \notag  \\
&&- \frac{1}{{\sigma _1 ^3 (B)}}\sum\limits_{i \in B}
{\sum\limits_{a,b = 1}^N {F^{ab} } [\sigma _1 (B)u_{iia}  - u_{ii}
\sum\limits_{j \in B} {u_{jja} } ]} [\sigma _1 (B)u_{iib}  - u_{ii}
\sum\limits_{j \in B}{u_{jjb} } ]\\
&&- \frac{1}{{\sigma _1 (B)}}\sum\limits_{\scriptstyle i,j \in B
\hfill \atop \scriptstyle i \ne j \hfill}  {\sum\limits_{a,b = 1}^N
{F^{ab} u_{ija} u_{ijb} }} \notag \\
&&+ O(\sum\limits_{i,j \in B} {\left| {\nabla u_{ij} } \right|}  +
\phi ). \notag
\end{eqnarray}
In fact, $u$ in (3.5) is $u_\varepsilon
(x)=u(x)+\frac{\varepsilon}{2} \left| {x'} \right|^2 $ defined as
above (we omit the subindex $\varepsilon$).
\end{theorem}
\textbf{Proof}. The proof is similar to the proof in [2]. We give
the main process.

Following the assumptions as above, and for a similar computation as
in [2], we have
\begin{align}
\sigma _1 (B) = O(\phi ),  u_{ii} = O(\phi ) \text{ for every } i
\in B.
\end{align}

Since $ \phi (x)= \sigma _{l + 1} (W) + q(W)$, then by the chain
rule we have
\begin{align}
 \sum\limits_{a,b = 1}^N {F^{ab} \phi _{ab} }  =&
\sum\limits_{a,b = 1}^N {F^{ab} [\sum\limits_{i,j} {\frac{{\partial
\phi }} {{\partial u_{ij} }}u_{ijab}  + \sum\limits_{i,j,k,l}
{\frac{{\partial ^2 \phi}} {{\partial u_{ij} \partial u_{kl} }}} u_{ija} u_{klb} ]} } \notag \\
 = &\sum\limits_{a,b = 1}^N {F^{ab} \sum\limits_{i,j} {[\frac{{\partial \sigma _{l + 1} (W)}}
{{\partial u_{ij} }} + \frac{{\partial q(W)}} {{\partial u_{ij}}}]u_{ijab} } } \\
& + \sum\limits_{a,b = 1}^N {F^{ab} \sum\limits_{i,j,k,l}
{[\frac{{\partial ^2 \sigma _{l + 1} (W)}} {{\partial u_{ij}
\partial u_{kl} }} + } \frac{{\partial ^2 q(W)}} {{\partial u_{ij}
\partial u_{kl} }}]u_{ija} u_{rsb} }. \notag
\end{align}

Since $W$ is diagonal and by lemma 2.2, the first term on the right
hand side of (3.7)
 is
\begin{align}
 \sum\limits_{a,b = 1}^N {F^{ab} \sum\limits_{i,j}
{\frac{{\partial \sigma _{l + 1} (W)}} {{\partial u_{ij} }}u_{ijab}
} }  = &\sum\limits_{a,b = 1}^N {F^{ab} \sum\limits_{i}
{\frac{{\partial \sigma _{l + 1} (W)}} {{\partial u_{ii} }}u_{iiab}
} } \notag  \\
= &\sum\limits_{a,b = 1}^N {F^{ab} \sum\limits_{i \in B} {\sigma _l
(G)u_{iiab} } }  + O(\phi).
\end{align}

Using Lemma 2.4 in [2], the second term on the right hand side of
(3.7) is
\begin{align}
\sum\limits_{a,b = 1}^N {F^{ab} \sum\limits_{i,j} {\frac{{\partial
q(W)}} {{\partial u_{ij} }}u_{ijab} } }  =& \sum\limits_{a,b = 1}^N
{F^{ab} \sum\limits_i {\frac{{\partial q(W)}} {{\partial u_{ii}
}}u_{iiab} } } \notag  \\
= &\sum\limits_{a,b = 1}^N {F^{ab} \sum\limits_{i \in B}
{\frac{{\sigma _1 ^2 (B\left| i \right.) - \sigma _2 (B\left| i
\right.)}} {{\sigma _1 ^2 (B)}}u_{iiab} } }+ O(\phi ).
\end{align}

As in [2], the third term on the right hand side of (3.7) is
\begin{align}
&\sum\limits_{a,b = 1}^N {F^{ab} \sum\limits_{i,j,k,l}
{\frac{{\partial ^2 \sigma _{l + 1} (W)}} {{\partial u_{ij} \partial
u_{kl} }}u_{ija} u_{klb} } } \notag \\
 = &\sum\limits_{a,b = 1}^N {F^{ab} [\sum\limits_{i \ne j} {\frac{{\partial ^2 \sigma _{l + 1}
(W)}} {{\partial u_{ii} \partial u_{jj} }}u_{iia} u_{jjb}  +
\sum\limits_{i \ne j} {\frac{{\partial ^2 \sigma _{l + 1} (W)}}
{{\partial u_{ij} \partial u_{ji} }}u_{ija} u_{jib} ]} } }\\
=  &- 2\sum\limits_{a,b = 1}^N {F^{ab} \sum\limits_{i \in B,j \in G}
{\sigma _{l - 1} (G\left| j \right.)u_{ija} u_{jib} } }   +
O(\sum\limits_{i,j \in B} {\left| {\nabla u_{ij} } \right|}  + \phi
). \notag
\end{align}

From Proposition 2.1 in [2], we can get
\begin{align}
&\sum\limits_{i,j,k,l} {\frac{{\partial ^2 q(W)}} {{\partial u_{ij}
\partial u_{kl} }}u_{ija} u_{klb} }  =  - 2\sum\limits_{i \in B,j
\in G} {\frac{{\sigma _1 ^2 (B\left| i \right.) - \sigma _2 (B\left|
i \right.)}} {{\sigma _1 ^2 (B)u_{jj} }}u_{ija} u_{ijb} } \notag \\
 &- \frac{1}{{\sigma _1 ^3 (B)}}\sum\limits_{i \in B} {[\sigma _1 (B)u_{iia}  -
u_{ii} \sum\limits_{j \in B} {u_{jja} } ][\sigma _1 (B)u_{iib}  -
u_{ii} \sum\limits_{j \in B} {u_{jjb} } ]}\\
&- \frac{1}{{\sigma _1 (B)}}\sum\limits_{\scriptstyle i,j \in B
\hfill \atop \scriptstyle i \ne j \hfill}  {u_{ija} u_{ijb} }  +
O(\sum\limits_{i,j \in B} {\left| {\nabla u_{ij} } \right|}  + \phi
).\notag
\end{align}
So the fourth term on the right hand side of (3.7) is
\begin{align}
 &\sum\limits_{a,b = 1}^N {F^{ab} \sum\limits_{i,j,k,l}
{\frac{{\partial ^2 q(W)}} {{\partial u_{ij} \partial u_{kl}
}}u_{ija} u_{klb} } } \notag \\
 =  &- 2\sum\limits_{a,b = 1}^N
{F^{ab} \sum\limits_{i \in B,j \in G} {\frac{{\sigma _1 ^2 (B\left|
i \right.) - \sigma _2 (B\left| i \right.)}} {{\sigma _1 ^2
(B)u_{jj}}}u_{ija} u_{ijb} } } \notag  \\
& - \frac{1} {{\sigma _1 ^3 (B)}}\sum\limits_{a,b = 1}^N {F^{ab}
\sum\limits_{i \in B} {[\sigma _1 (B)u_{iia}  - u_{ii}
\sum\limits_{j \in B} {u_{jja} } ][\sigma _1 (B)u_{iib}  - u_{ii}
\sum\limits_{j \in B} {u_{jjb} } ]} }  \\
& - \frac{1} {{\sigma _1 (B)}}\sum\limits_{a,b = 1}^N {F^{ab}
\sum\limits_{\scriptstyle i,j \in B \hfill \atop \scriptstyle i \ne
j \hfill}  {u_{ija} u_{ijb} } }  + O(\sum\limits_{i,j \in B} {\left|
{\nabla u_{ij} } \right|}  + \phi ). \notag
\end{align}
Substitute (3.8), (3.9), (3.10) and (3.12) into (3.7), then we
obtain(3.5).

\subsection{calculation on structure condition}
Now we discuss the structure condition (1.4). We write
$ A= \left(
{\begin{matrix}
   a & b  \\
   {b^T } & c  \\
 \end{matrix} } \right)
$ and $a^{-1}=(a^{ij})$, where $a=(a_{ij}) \in \mathcal{S}^{N'}$,
$b=(b_{k \alpha}) \in \mathbb{R}^{N' \times N''}$ and $c=(c_{\alpha
\beta}) \in \mathcal{S}^{N''}$.
\begin{lemma}
The condition (1.4) is equivalent to
\begin{align}
&\sum\limits_{a,b,c,d = 1}^N {F^{ab,cd} (A,p,u,x)X_{ab} X_{cd} }  +
2\sum\limits_{a,b = 1}^N {\sum\limits_{k,l = 1}^{N'} {F^{ab} a^{kl}
X_{ka} X_{lb} } }  + 2\sum\limits_{a,b = 1}^N {\sum\limits_{\alpha =
N' + 1}^N {F^{ab,p_\alpha  } X_{ab} X_\alpha  } }  \\
 &+ 2\sum\limits_{a,b = 1}^N {F^{ab,u} X_{ab} Y}  + 2\sum\limits_{a,b = 1}^N
 {\sum\limits_{i = 1}^{N'} {F^{ab,x_i } X_{ab} Z_i } }  + \sum\limits_{\alpha ,\beta  = N' + 1}^N
 {F^{p_\alpha  ,p_\beta  } X_\alpha  X_\beta  }  + 2\sum\limits_{\alpha  = N' + 1}^N
 {F^{p_\alpha  ,u} X_\alpha Y} \notag  \\
&+ 2\sum\limits_{\alpha  = N' + 1}^N {\sum\limits_{i = 1}^{N'}
{F^{p_\alpha  ,x_i } X_\alpha  Z_i } }
 + F^{u,u} Y^2  + 2\sum\limits_{i = 1}^{N'} {F^{u,x_i } YZ_i }  + \sum\limits_{i,j = 1}^{N'}
 {F^{x_i ,x_j } Z_i Z_j }  \geqslant 0, \notag
\end{align}
for every $ \widetilde X = ((X_{ab} ),(X_\alpha  ),Y,(Z_i )) \in
\mathcal{S}^N \times \mathbb{R}^{N''}  \times \mathbb{R} \times
\mathbb{R}^{N'}$.
\end{lemma}
Proof. We denote $G(a,b,c,p'',u,x')=F(\left( {\begin{matrix}
   a^{-1} & a^{-1}b  \\
   {({a^{-1}b})^T } & c+ b^Ta^{-1}b \\
 \end{matrix} } \right),p',p'',u,x',x'') $, and we let $ 1 \leqslant
 i,j,k,l,m,n,s,t \leqslant N' $, $ N' \leqslant \alpha ,\beta ,\gamma ,\eta ,\xi
,\zeta \leqslant N'' $, and $ 1 \leqslant a,b,c,d \leqslant N $.
Then condition (1.4) is equivalent to
\begin{align}
&\sum\limits_{i,j,k,l} {\frac{{\partial ^2 G}} {{\partial a_{ij}
\partial a_{kl} }}X_{ij} X_{kl} }  + 2\sum\limits_{i,j,k,\alpha }
{\frac{{\partial ^2 G}} {{\partial a_{ij} \partial b_{k\alpha }
}}X_{ij} Y_{k\alpha } }  + 2\sum\limits_{i,j,\alpha ,\beta }
{\frac{{\partial ^2 G}} {{\partial a_{ij} \partial c_{\alpha \beta }
}}X_{ij} Z_{\alpha \beta } }  + 2\sum\limits_{i,j,\alpha }
{\frac{{\partial ^2 G}} {{\partial a_{ij} \partial p_\alpha }}X_{ij}
X_\alpha  } \\
& + 2\sum\limits_{i,j} {\frac{{\partial ^2 G}} {{\partial a_{ij}
\partial u}}X_{ij} Y}  + 2\sum\limits_{i,j,k} {\frac{{\partial ^2
G}} {{\partial a_{ij} \partial x_k }}X_{ij} Z_k }  +
\sum\limits_{k,\alpha ,l,\beta } {\frac{{\partial ^2 G}} {{\partial
b_{k\alpha } \partial b_{l\beta } }}Y_{k\alpha } Y_{l\beta } }  +
2\sum\limits_{k,\alpha ,\gamma ,\eta } {\frac{{\partial ^2 G}}
{{\partial b_{k\alpha } \partial c_{\gamma \eta } }}Y_{k\alpha }
Z_{\gamma \eta } } \notag  \\
&+ 2\sum\limits_{k,\alpha ,\beta } {\frac{{\partial ^2 G}}
{{\partial b_{k\alpha } \partial p_\beta  }}Y_{k\alpha } X_\beta  }
+ 2\sum\limits_{k,\alpha } {\frac{{\partial ^2 G}} {{\partial
b_{k\alpha } \partial u}}Y_{k\alpha } Y}  + 2\sum\limits_{k,\alpha
,l} {\frac{{\partial ^2 G}} {{\partial b_{k\alpha } \partial x_l
}}Y_{k\alpha } Z_l }  + \sum\limits_{\alpha ,\beta ,\gamma ,\eta }
{\frac{{\partial ^2 G}} {{\partial c_{\alpha \beta } \partial
c_{\gamma \eta } }}Z_{\alpha \beta } Z_{\gamma \eta } } \notag  \\
 &+ 2\sum\limits_{\alpha ,\beta ,\gamma } {\frac{{\partial ^2 G}}
{{\partial c_{\alpha \beta } \partial p_\gamma  }}Z_{\alpha \beta }
X_\gamma  }  + 2\sum\limits_{\alpha ,\beta} {\frac{{\partial ^2 G}}
{{\partial c_{\alpha \beta } \partial u}}Z_{\alpha \beta } Y}  +
2\sum\limits_{\alpha ,\beta ,l} {\frac{{\partial ^2 G}} {{\partial
c_{\alpha \beta } \partial x_l }}Z_{\alpha \beta } Z_l  + }
\sum\limits_{\alpha ,\gamma } {\frac{{\partial ^2 G}} {{\partial
p_\alpha  \partial p_\gamma
}}X_\alpha  X_\gamma  } \notag  \\
& + 2\sum\limits_\alpha  {\frac{{\partial ^2 G}} {{\partial p_\alpha
\partial u}}X_\alpha  Y}  + 2\sum\limits_{\alpha,l}{\frac{{\partial ^2
G}} {{\partial p_\alpha  \partial x_l }}X_\alpha Z_l }  +
\frac{{\partial ^2 G}} {{\partial u\partial u}}Y^2  + 2\sum\limits_l
{\frac{{\partial ^2 G}} {{\partial u\partial x_l }}YZ_l }  +
\sum\limits_{k,l} {\frac{{\partial ^2 G}} {{\partial x_k
\partial x_l }}Z_k Z_l }  \geqslant 0, \notag
\end{align}
for every $ ((X_{ij}), (Y_{k\alpha }), (Z_{\alpha \beta }),
(X_\alpha), Y, (Z_i) ) \in \mathcal{S}^{N'}  \times \mathbb{R}^{N'
\times N''} \times \mathcal{S}^{N''} \times \mathbb{R}^{N''} \times
\mathbb{R} \times \mathbb{R}^{N'} $.

To get the equivalent condition (3.13), we shall represent all the
derivatives of $G$ in (3.14) by the derivatives of $F$.

Suppose $ a^{ - 1} b = (B_{k\alpha } ) = (\sum\limits_l {a^{kl}
b_{l\alpha } } ) $, and $ c + b^T a^{ - 1} b = (C_{\alpha \beta } )
= (c_{\alpha \beta } + \sum\limits_{k,l} {b_{k\alpha } a^{kl}
b_{l\beta } } ) $, then $G(a,b,c,p'',u,x')=F(\left( {\begin{matrix}
   a^{-1} & (B_{k\alpha} )  \\   {{(B_{k\alpha}) }^T } & (C_{\alpha \beta } ) \\
 \end{matrix} } \right),p',p'',u,x',x'') $.\\

A direct computation yields
\begin{align}
&\frac{{\partial G}} {{\partial a_{ij} }} = \sum\limits_{k,l}
{F^{kl} \frac{{\partial a^{kl} }} {{\partial a_{ij} }} + }
\sum\limits_{k,\beta } {F^{k\beta } \frac{{\partial B_{k\beta } }}
{{\partial a_{ij} }}}  + \sum\limits_{\alpha ,l} {F^{\alpha l}
\frac{{\partial B_{l\alpha } }} {{\partial a_{ij} }}}  +
\sum\limits_{\alpha ,\beta } {F^{\alpha \beta } \frac{{\partial
C_{\alpha \beta } }} {{\partial a_{ij} }}}, \\
&\frac{{\partial G}} {{\partial b_{k\beta } }} =
\sum\limits_{m,\alpha } {F^{m\alpha } \frac{{\partial B_{m\alpha }
}} {{\partial b_{k\beta } }}}  + \sum\limits_{\alpha ,n} {F^{\alpha
n} \frac{{\partial B_{n\alpha } }} {{\partial b_{k\beta } }}}  +
\sum\limits_{\gamma ,\eta } {F^{\gamma \eta } \frac{{\partial
C_{\gamma \eta } }} {{\partial b_{k\beta } }}}.
\end{align}
So we have the second derivatives of $G$ in (3.14) as follows. The
derivatives of $G$ in the last ten terms are simple,
\begin{align*}
& \frac{{\partial ^2 G}} {{\partial x_k \partial x_l }} = F^{x_k
,x_l }, \frac{{\partial ^2 G}} {{\partial u\partial x_l }} =
F^{u,x_l}, \frac{{\partial ^2 G}} {{\partial u\partial u}} = F^{u,u},\\
& \frac{{\partial ^2 G}} {{\partial p_\alpha  \partial x_l }} =
F^{p_\alpha  ,x_l }, \frac{{\partial ^2 G}} {{\partial p_\alpha
\partial u}} = F^{p_\alpha  ,u}, \frac{{\partial ^2 G}} {{\partial
p_\alpha  \partial p_\gamma  }} = F^{p_\alpha  ,p_\gamma  },\\
& \frac{{\partial ^2 G}} {{\partial c_{\alpha \beta } \partial x_l
}} = F^{\alpha \beta ,x_l } ,\frac{{\partial ^2 G}} {{\partial
c_{\alpha \beta } \partial u}} = F^{\alpha \beta ,u},\\
& \frac{{\partial ^2 G}} {{\partial c_{\alpha \beta } \partial
p_\gamma  }} = F^{\alpha \beta ,p_\gamma  } ,\frac{{\partial ^2 G}}
{{\partial c_{\alpha \beta } \partial c_{\gamma \eta } }} =
F^{\alpha \beta ,\gamma \eta }.
\end{align*}
From (3.15), we can get the derivatives of $G$ in the third-sixth
terms of (3.14)
$$
\frac{{\partial ^2 G}} {{\partial a_{ij} \partial c_{\gamma \eta }
}} = \sum\limits_{k,l} {F^{kl,\gamma \eta } \frac{{\partial a^{kl}
}} {{\partial a_{ij} }} + } \sum\limits_{k,\beta } {F^{k\beta
,\gamma \eta } \frac{{\partial B_{k\beta } }} {{\partial a_{ij} }}}
+ \sum\limits_{\alpha ,l} {F^{\alpha l,\gamma \eta } \frac{{\partial
B_{l\alpha } }} {{\partial a_{ij} }}}  + \sum\limits_{\alpha ,\beta
} {F^{\alpha \beta ,\gamma \eta } \frac{{\partial C_{\alpha \beta }
}} {{\partial a_{ij} }}},
$$

$$
\frac{{\partial ^2 G}} {{\partial a_{ij} \partial p_\gamma  }} =
\sum\limits_{k,l} {F^{kl,p_\gamma  } \frac{{\partial a^{kl} }}
{{\partial a_{ij} }} + } \sum\limits_{k,\beta } {F^{k\beta ,p_\gamma
} \frac{{\partial B_{k\beta } }} {{\partial a_{ij} }}}  +
\sum\limits_{\alpha ,l} {F^{\alpha l,p_\gamma  } \frac{{\partial
B_{l\alpha } }} {{\partial a_{ij} }}}  + \sum\limits_{\alpha ,\beta
} {F^{\alpha \beta ,p_\gamma  } \frac{{\partial C_{\alpha \beta } }}
{{\partial a_{ij} }}},
$$

$$
\frac{{\partial ^2 G}} {{\partial a_{ij} \partial u}} =
\sum\limits_{k,l} {F^{kl,u} \frac{{\partial a^{kl} }} {{\partial
a_{ij} }} + } \sum\limits_{k,\beta } {F^{k\beta ,u} \frac{{\partial
B_{k\beta } }} {{\partial a_{ij} }}}  + \sum\limits_{\alpha ,l}
{F^{\alpha l,u} \frac{{\partial B_{l\alpha } }} {{\partial a_{ij}
}}}  + \sum\limits_{\alpha ,\beta } {F^{\alpha \beta ,u}
\frac{{\partial C_{\alpha \beta } }} {{\partial a_{ij} }}},
$$

$$
\frac{{\partial ^2 G}} {{\partial a_{ij} \partial x_m }} =
\sum\limits_{k,l} {F^{kl,x_m } \frac{{\partial a^{kl} }} {{\partial
a_{ij} }} + } \sum\limits_{k,\beta } {F^{k\beta ,x_m }
\frac{{\partial B_{k\beta } }} {{\partial a_{ij} }}}  +
\sum\limits_{\alpha ,l} {F^{\alpha l,x_m } \frac{{\partial
B_{l\alpha } }} {{\partial a_{ij} }}}  + \sum\limits_{\alpha ,\beta
} {F^{\alpha \beta ,x_m } \frac{{\partial C_{\alpha \beta } }}
{{\partial a_{ij} }}}.
$$
From (3.16), we can get the derivatives of $G$ in the
eighth-eleventh terms of (3.14)
$$
\frac{{\partial ^2 G}} {{\partial b_{k\beta } \partial c_{\xi \zeta
} }} = \sum\limits_{m,\alpha } {F^{m\alpha ,\xi \zeta }
\frac{{\partial B_{m\alpha } }} {{\partial b_{k\beta } }}}  +
\sum\limits_{\alpha ,n} {F^{\alpha n,\xi \zeta } \frac{{\partial
B_{n\alpha } }} {{\partial b_{k\beta } }}}  + \sum\limits_{\gamma
,\eta } {F^{\gamma \eta ,\xi \zeta } \frac{{\partial C_{\gamma \eta
} }} {{\partial b_{k\beta } }}},
$$

$$
\frac{{\partial ^2 G}} {{\partial b_{k\beta } \partial p_\zeta  }} =
\sum\limits_{m,\alpha } {F^{m\alpha ,p_\zeta  } \frac{{\partial
B_{m\alpha } }} {{\partial b_{k\beta } }}}  + \sum\limits_{\alpha
,n} {F^{\alpha n,p_\zeta  } \frac{{\partial B_{n\alpha } }}
{{\partial b_{k\beta } }}}  + \sum\limits_{\gamma ,\eta } {F^{\gamma
\eta ,p_\zeta  } \frac{{\partial C_{\gamma \eta } }} {{\partial
b_{k\beta } }}},
$$

$$
\frac{{\partial ^2 G}} {{\partial b_{k\beta } \partial u}} =
\sum\limits_{m,\alpha } {F^{m\alpha ,u} \frac{{\partial B_{m\alpha }
}} {{\partial b_{k\beta } }}}  + \sum\limits_{\alpha ,n} {F^{\alpha
n,u} \frac{{\partial B_{n\alpha } }} {{\partial b_{k\beta } }}}  +
\sum\limits_{\gamma ,\eta } {F^{\gamma \eta ,u} \frac{{\partial
C_{\gamma \eta } }} {{\partial b_{k\beta } }}},
$$

$$
\frac{{\partial ^2 G}} {{\partial b_{k\beta } \partial x_i }} =
\sum\limits_{m,\alpha } {F^{m\alpha ,x_i } \frac{{\partial
B_{m\alpha } }} {{\partial b_{k\beta } }}}  + \sum\limits_{\alpha
,n} {F^{\alpha n,x_i } \frac{{\partial B_{n\alpha } }} {{\partial
b_{k\beta } }}}  + \sum\limits_{\gamma ,\eta } {F^{\gamma \eta ,x_i
} \frac{{\partial C_{\gamma \eta } }} {{\partial b_{k\beta } }}}.
$$
Also from (3.15) we can get the derivative of $G$ in the first term
of (3.14)
\begin{align*}
\frac{{\partial ^2 G}} {{\partial a_{ij} \partial a_{mn} }}& =
\sum\limits_{k,l} {F^{kl} \frac{{\partial ^2 a^{kl} }} {{\partial
a_{ij} \partial a_{mn} }} + } \sum\limits_{k,\beta } {F^{k\beta }
\frac{{\partial ^2 B_{k\beta } }} {{\partial a_{ij} \partial a_{mn}
}}}  + \sum\limits_{\alpha ,l} {F^{\alpha l} \frac{{\partial ^2
B_{l\alpha } }} {{\partial a_{ij} \partial a_{mn} }}}  +
\sum\limits_{\alpha ,\beta } {F^{\alpha \beta } \frac{{\partial ^2
C_{\alpha \beta } }} {{\partial a_{ij} \partial a_{mn} }}}\\
 &+ \sum\limits_{k,l} {\frac{{\partial a^{kl} }}
{{\partial a_{ij} }}[\sum\limits_{s,t} {F^{kl,st} \frac{{\partial
a^{st} }} {{\partial a_{mn} }} + } \sum\limits_{s,\eta }
{F^{kl,s\eta } \frac{{\partial B_{s\eta } }} {{\partial a_{mn} }}} +
\sum\limits_{\gamma ,t} {F^{kl,\gamma t} \frac{{\partial B_{t\gamma
} }} {{\partial a_{mn} }}}  + \sum\limits_{\gamma, \eta }
{F^{kl,\gamma \eta } \frac{{\partial C_{\gamma \eta } }} {{\partial
a_{mn} }}} } ]\\
 &+ \sum\limits_{k,\beta } {\frac{{\partial B_{k\beta } }}
{{\partial a_{ij} }}[\sum\limits_{s,t} {F^{k\beta ,st}
\frac{{\partial a^{st} }} {{\partial a_{mn} }} + }
\sum\limits_{s,\eta } {F^{k\beta ,s\eta } \frac{{\partial B_{s\eta }
}} {{\partial a_{mn} }}}  + \sum\limits_{\gamma ,t} {F^{k\beta
,\gamma t} \frac{{\partial B_{t\gamma } }} {{\partial a_{mn} }}}  +
\sum\limits_{\gamma,\eta } {F^{k\beta ,\gamma \eta } \frac{{\partial
C_{\gamma \eta } }} {{\partial a_{mn} }}} ]}\\
& + \sum\limits_{\alpha ,l} {\frac{{\partial B_{l\alpha } }}
{{\partial a_{ij} }}[\sum\limits_{s,t} {F^{\alpha l,st}
\frac{{\partial a^{st} }} {{\partial a_{mn} }} + }
\sum\limits_{s,\eta } {F^{\alpha l,s\eta } \frac{{\partial B_{s\eta
} }} {{\partial a_{mn} }}}  + \sum\limits_{\gamma ,t} {F^{\alpha
l,\gamma t} \frac{{\partial B_{t\gamma } }} {{\partial a_{mn} }}}  +
\sum\limits_{\gamma, \eta } {F^{\alpha l,\gamma \eta }
\frac{{\partial C_{\gamma \eta } }} {{\partial a_{mn} }}} ]}\\
 &+ \sum\limits_{\alpha ,\beta } {\frac{{\partial C_{\alpha \beta } }}
{{\partial a_{ij} }}[\sum\limits_{s,t} {F^{\alpha \beta ,st}
\frac{{\partial a^{st} }} {{\partial a_{mn} }} + }
\sum\limits_{s,\eta } {F^{\alpha \beta ,s\eta } \frac{{\partial
B_{s\eta } }} {{\partial a_{mn} }}}  + \sum\limits_{\gamma ,t}
{F^{\alpha \beta ,\gamma t} \frac{{\partial B_{t\gamma } }}
{{\partial a_{mn} }}}  + \sum\limits_{\gamma, \eta } {F^{\alpha
\beta ,\gamma \eta } \frac{{\partial C_{\gamma \eta } }} {{\partial
a_{mn} }}} ]},
\end{align*}
and the derivative of $G$ in the second term in (3.14)
\begin{align*}
\frac{{\partial ^2 G}} {{\partial a_{ij} \partial b_{m\eta } }} &=
\sum\limits_{k,\beta } {F^{k\beta } \frac{{\partial ^2 B_{k\beta }
}} {{\partial a_{ij} \partial b_{m\eta } }}}  + \sum\limits_{\alpha
,l} {F^{\alpha l} \frac{{\partial ^2 B_{l\alpha } }} {{\partial
a_{ij} \partial b_{m\eta } }}}  + \sum\limits_{\alpha ,\beta }
{F^{\alpha \beta } \frac{{\partial ^2 C_{\alpha \beta } }}
{{\partial a_{ij} \partial b_{m\eta } }}}\\
& + \sum\limits_{k,l} {\frac{{\partial a^{kl} }} {{\partial a_{ij}
}}[\sum\limits_{s,\zeta } {F^{kl,s\zeta } \frac{{\partial B_{s\zeta
} }} {{\partial b_{m\eta } }}}  + \sum\limits_{\xi ,t} {F^{kl,\xi t}
\frac{{\partial B_{t\xi } }} {{\partial b_{m\eta } }}}  +
\sum\limits_{\xi ,\zeta } {F^{kl,\xi \zeta } \frac{{\partial C_{\xi
\zeta } }} {{\partial b_{m\eta } }}}
} ]\\
& + \sum\limits_{k,\beta } {\frac{{\partial B_{k\beta } }}
{{\partial a_{ij} }}[\sum\limits_{s,\zeta } {F^{k\beta ,s\zeta }
\frac{{\partial B_{s\zeta } }} {{\partial b_{m\eta } }}}  +
\sum\limits_{\xi ,t} {F^{k\beta ,\xi t} \frac{{\partial B_{t\xi } }}
{{\partial b_{m\eta } }}}  + \sum\limits_{\xi ,\zeta } {F^{k\beta
,\xi \zeta } \frac{{\partial C_{\xi \zeta } }} {{\partial b_{m\eta }
}}} ]}\\
& + \sum\limits_{\alpha ,l} {\frac{{\partial B_{l\alpha } }}
{{\partial a_{ij} }}[\sum\limits_{s,\zeta } {F^{\alpha l,s\zeta }
\frac{{\partial B_{s\zeta } }} {{\partial b_{m\eta } }}}  +
\sum\limits_{\xi ,t} {F^{\alpha l,\xi t} \frac{{\partial B_{t\xi }
}} {{\partial b_{m\eta } }}}  + \sum\limits_{\xi ,\zeta } {F^{\alpha
l,\xi \zeta } \frac{{\partial C_{\xi \zeta } }} {{\partial b_{m\eta
} }}} ]}\\
& + \sum\limits_{\alpha ,\beta } {\frac{{\partial C_{\alpha \beta }
}} {{\partial a_{ij} }}[\sum\limits_{s,\zeta } {F^{\alpha \beta
,s\zeta } \frac{{\partial B_{s\zeta } }} {{\partial b_{m\eta } }}} +
\sum\limits_{\xi ,t} {F^{\alpha \beta ,\xi t} \frac{{\partial
B_{t\xi } }} {{\partial b_{m\eta } }}}  + \sum\limits_{\xi ,\zeta }
{F^{\alpha \beta ,\xi \zeta } \frac{{\partial C_{\xi \zeta } }}
{{\partial b_{m\eta } }}} ]}.
\end{align*}
From (3.16), we can get the derivative of $G$ in the seventh term of
(3.14)
\begin{align*}
\frac{{\partial ^2 G}} {{\partial b_{k\beta } \partial b_{l\alpha }
}} &= \sum\limits_{\gamma ,\eta } {F^{\gamma \eta } \frac{{\partial
^2 C_{\gamma \eta } }} {{\partial b_{k\beta } \partial b_{l\alpha }
}}}\\
 &+ \sum\limits_{m,\eta } {\frac{{\partial B_{m\eta } }}
{{\partial b_{k\beta } }}} [\sum\limits_{s,\zeta } {F^{m\eta ,s\zeta
} \frac{{\partial B_{s\zeta } }} {{\partial b_{l\alpha } }}}  +
\sum\limits_{\xi ,t} {F^{m\eta ,\xi t} \frac{{\partial B_{t\xi } }}
{{\partial b_{l\alpha } }}}  + \sum\limits_{\xi ,\zeta } {F^{m\eta
,\xi \zeta } \frac{{\partial C_{\xi \zeta } }} {{\partial b_{l\alpha
} }}} ]\\
& + \sum\limits_{\gamma ,n} {\frac{{\partial B_{n\gamma } }}
{{\partial b_{k\beta } }}} [\sum\limits_{s,\zeta } {F^{\gamma
n,s\zeta } \frac{{\partial B_{s\zeta } }} {{\partial b_{l\alpha }
}}} + \sum\limits_{\xi ,t} {F^{\gamma n,\xi t} \frac{{\partial
B_{t\xi } }} {{\partial b_{l\alpha } }}}  + \sum\limits_{\xi ,\zeta
} {F^{\gamma n,\xi \zeta } \frac{{\partial C_{\xi \zeta } }}
{{\partial b_{l\alpha } }}} ]\\
& + \sum\limits_{\gamma ,\eta } {\frac{{\partial C_{\gamma \eta } }}
{{\partial b_{k\beta } }}} [\sum\limits_{s,\zeta } {F^{\gamma \eta
,s\zeta } \frac{{\partial B_{s\zeta } }} {{\partial b_{l\alpha } }}}
+ \sum\limits_{\xi ,t} {F^{\gamma \eta ,\xi t} \frac{{\partial
B_{t\xi } }} {{\partial b_{l\alpha } }}}  + \sum\limits_{\xi ,\zeta
} {F^{\gamma \eta ,\xi \zeta } \frac{{\partial C_{\xi \zeta } }}
{{\partial b_{l\alpha } }}} ].
\end{align*}
So we denote
\begin{align}
&\widetilde X_{kl}  = \sum\limits_{i,j} {\frac{{\partial a^{kl} }}
{{\partial a_{ij} }}X_{ij} } ,\widetilde X_{\beta k}  = \widetilde
X_{k\beta }  = \sum\limits_{i,j} {\frac{{\partial B_{k\beta } }}
{{\partial a_{ij} }}X_{ij} } ,\widetilde X_{\alpha \beta }  =
\sum\limits_{i,j} {\frac{{\partial C_{\alpha \beta } }} {{\partial
a_{ij} }}X_{ij} },\\
& \widetilde Y_{kl}  = 0,\widetilde Y_{\beta k}  = \widetilde
Y_{k\beta } = \sum\limits_{m,\eta } {\frac{{\partial B_{k\beta } }}
{{\partial b_{m\eta } }}Y_{m\eta } } ,\widetilde Y_{\alpha \beta } =
\sum\limits_{m,\eta } {\frac{{\partial C_{\alpha \beta } }}
{{\partial b_{m\eta } }}Y_{m\eta } },\\
& \widetilde Z_{kl}  = 0,\widetilde Z_{\beta k}  = \widetilde
Z_{k\beta } = 0,\widetilde Z_{\alpha \beta }  = Z_{\alpha \beta }.
\end{align}

From the above calculation, and (3.17)-(3.19), we can get the first
term of (3.14)
\begin{align}
\sum\limits_{i,j,m,n} &{\frac{{\partial ^2 G}} {{\partial a_{ij}
\partial a_{mn} }} X_{ij} X_{mn} } \notag \\
 =&  \sum\limits_{k,l} {F^{kl}
\sum\limits_{i,j,m,n} {\frac{{\partial ^2 a^{kl} }} {{\partial
a_{ij} \partial a_{mn} }}} X_{ij} X_{mn}  + } \sum\limits_{k,\beta }
{F^{k\beta } \sum\limits_{i,j,m,n} {\frac{{\partial ^2 B_{k\beta }
}} {{\partial a_{ij} \partial a_{mn} }}} } X_{ij} X_{mn} \notag \\
 &+ \sum\limits_{\alpha ,l} {F^{\alpha l} \sum\limits_{i,j,m,n} {\frac{{\partial ^2 B_{l\alpha } }}
{{\partial a_{ij} \partial a_{mn} }}X_{ij} X_{mn} } }  +
\sum\limits_{\alpha ,\beta } {F^{\alpha \beta }
\sum\limits_{i,j,m,n} {\frac{{\partial ^2 C_{\alpha \beta } }}
{{\partial a_{ij} \partial a_{mn} }}} X_{ij} X_{mn} } \notag \\
& + \sum\limits_{a,b,c,d = 1}^N {F^{ab,cd} \widetilde X_{ab}
\widetilde X_{cd} } \notag \\
=& 2\sum\limits_{a,b} {F^{ab} } \sum\limits_{ij} {a_{ij} }
\widetilde X_{ia} \widetilde X_{jb}+ \sum\limits_{a,b,c,d = 1}^N
{F^{ab,cd} \widetilde X_{ab} \widetilde X_{cd} },
\end{align}

the second term of (3.14)
\begin{align}
\sum\limits_{i,j,m,\eta } &{\frac{{\partial ^2 G}} {{\partial a_{ij}
\partial b_{m\eta } }}X_{ij} Y_{m\eta } } \notag \\
=& \sum\limits_{k,\beta }
{F^{k\beta } \sum\limits_{i,j,m,\eta } {\frac{{\partial ^2 B_{k\beta
} }} {{\partial a_{ij} \partial b_{m\eta } }}X_{ij} Y_{m\eta } } } +
\sum\limits_{\alpha ,l} {F^{\alpha l} \sum\limits_{i,j,m,\eta }
{\frac{{\partial ^2 B_{l\alpha } }} {{\partial a_{ij} \partial
b_{m\eta } }}X_{ij} Y_{m\eta } } } \notag \\
&+ \sum\limits_{\alpha ,\beta } {F^{\alpha \beta }
\sum\limits_{i,j,m,\eta } {\frac{{\partial ^2 C_{\alpha \beta } }}
{{\partial a_{ij} \partial b_{m\eta } }}X_{ij} Y_{m\eta } } }  +
\sum\limits_{a,b,c,d = 1}^N {F^{ab,cd} \widetilde X_{ab} \widetilde Y_{cd} } \notag \\
=& \sum\limits_{k,\beta } {F^{k\beta } \sum\limits_{i,j} {a_{ij}\widetilde X_{ik}
\widetilde Y_{j\beta } } }+ \sum\limits_{\alpha ,l} {F^{\alpha l} \sum\limits_{i,j}
{a_{ij} \widetilde X_{il } \widetilde Y_{j\alpha} } }\notag \\
&+ \sum\limits_{\alpha ,\beta } {F^{\alpha \beta } \sum\limits_{i,j}
{a_{ij} (\widetilde X_{i\alpha } \widetilde Y_{j\beta } + \widetilde
X_{i\beta } \widetilde Y_{j\alpha } )} }  + \sum\limits_{a,b,c,d =
1}^N {F^{ab,cd}\widetilde X_{ab} \widetilde Y_{cd} }\notag \\
= &\sum\limits_{a,b} {F^{ab} \sum\limits_{i,j} {a_{ij} (\widetilde
X_{ia} \widetilde Y_{jb} + \widetilde X_{ib} \widetilde Y_{ja} )} }
+ \sum\limits_{a,b,c,d = 1}^N {F^{ab,cd} \widetilde X_{ab}
\widetilde Y_{cd} },
\end{align}
and the seventh term of (3.14)
\begin{align}
\sum\limits_{k,\beta ,l,\alpha } {\frac{{\partial ^2 G}} {{\partial
b_{k\beta } \partial b_{l\alpha } }}Y_{k\beta } Y_{l\alpha } }  =
&\sum\limits_{\gamma ,\eta } {F^{\gamma \eta } \sum\limits_{k,\beta
,l,\alpha } {\frac{{\partial ^2 C_{\gamma \eta } }} {{\partial
b_{k\beta } \partial b_{l\alpha } }}Y_{k\beta } Y_{l\alpha } } }  +
\sum\limits_{a,b,c,d = 1}^N {F^{ab,cd} \widetilde Y_{ab}
\widetilde Y_{cd} } \notag \\
=& 2\sum\limits_{\gamma ,\eta } {F^{\gamma \eta } \sum\limits_{i,j}
{a_{ij} \widetilde Y_{i\gamma } \widetilde Y_{j\eta}}}
+\sum\limits_{a,b,c,d = 1}^N {F^{ab,cd} \widetilde Y_{ab} \widetilde
Y_{cd}} \notag \\
= &2\sum\limits_{a,b} {F^{ab} \sum\limits_{i,j} {a_{ij} \widetilde
Y_{ia} \widetilde Y_{jb} } } + \sum\limits_{a,b,c,d = 1}^N
{F^{ab,cd} \widetilde Y_{ab} \widetilde Y_{cd} }.
\end{align}
Also we obtain the third-sixth terms in (3.14)
\begin{align}
&\sum\limits_{i,j,\gamma ,\eta } {\frac{{\partial ^2 G}} {{\partial
a_{ij} \partial c_{\gamma \eta } }}X_{ij} Z_{\gamma \eta } }  =
\sum\limits_{a,b,c,d = 1}^N {F^{ab,cd} \widetilde X_{ab} \widetilde
Z_{cd} },\\
& \sum\limits_{i,j,\gamma } {\frac{{\partial ^2 G}} {{\partial
a_{ij}\partial p_\gamma  }}X_{ij} X_\gamma  }  = \sum\limits_{a,b =
1}^N {\sum\limits_\gamma  {F^{ab,p_\gamma  } \widetilde X_{ab}
X_\gamma }},\\
&\sum\limits_{i,j} {\frac{{\partial ^2 G}} {{\partial a_{ij}
\partial u}}X_{ij} Y}  = \sum\limits_{a,b = 1}^N {F^{ab,u}
\widetilde X_{ab} Y}, \\
&\sum\limits_{i,j,k}{\frac{{\partial ^2 G}} {{\partial a_{ij}
\partial x_k }}X_{ij} Z_k }  = \sum\limits_{a,b = 1}^N
{\sum\limits_k {F^{ab,x_k } \widetilde X_{ab} Z_k } },
\end{align}
and the eighth-eleventh terms in (3.14)
\begin{align}
& \sum\limits_{k,\beta ,\xi ,\zeta } {\frac{{\partial ^2 G}}
{{\partial b_{k\beta } \partial c_{\xi \zeta } }}Y_{k\beta } Z_{\xi
\zeta } }  = \sum\limits_{a,b,c,d = 1}^N {F^{ab,cd} \widetilde
Y_{ab} \widetilde Z_{cd} },\\
&\sum\limits_{k,\beta ,\xi } {\frac{{\partial ^2 G}} {{\partial
b_{k\beta } \partial p_\zeta  }}Y_{k\beta } X_\zeta  }  =
\sum\limits_{a,b = 1}^N {\sum\limits_\zeta  {F^{ab,p_\zeta  }
\widetilde Y_{ab} X_\zeta  } },\\
& \sum\limits_{k,\beta } {\frac{{\partial ^2 G}} {{\partial
b_{k\beta } \partial u}}Y_{k\beta } Y}  = \sum\limits_{a,b = 1}^N
{F^{ab,u} \widetilde Y_{ab} Y},\\
& \sum\limits_{k,\beta ,i} {\frac{{\partial ^2 G}} {{\partial
b_{k\beta } \partial x_i }}Y_{k\beta } Z_i }  = \sum\limits_{a,b =
1}^N {\sum\limits_i {F^{ab,x_i } \widetilde Y_{ab} Z_i } }.
\end{align}
So let$ \widetilde X = ((\widetilde X_{ab}  + \widetilde Y_{ab}  +
\widetilde Z_{ab} ),(X_\alpha  ),Y,(Z_i ))$, then we can obtain
(3.13). Also the equivalence holds.

\section{ structure condition and the proof of theorem 1.2 }
In this section, we prove Theorem 1.2 using a strong maximum
principle and Lemma 3.3. Also Corollary 1.4 holds directly from the
proof.

We denote $\mathcal{S}^n$ to be the set of all real symmetric $n
\times n$ matrices, and denote $\mathcal{S}_+^n \in \mathcal{S}^n$
to be the set of all positive definite symmetric $n \times n$
matrices. Let $\mathbb{O}_n$ be the space consisting all $n \times
n$ orthogonal matrices and $I_{N''}$ be the $N'' \times N''$
identity matrix. We define

$$
\mathcal{S}_{N' - 1}  = \{ \left( {\begin{matrix}
   {Q\left( {\begin{matrix}
   0 & 0  \\
   0 & B  \\
 \end{matrix} } \right)Q^T } & {Qb}  \\
   {b^TQ^T } & c  \\
 \end{matrix} } \right) \in \mathcal{S}^N \left| {\forall b \in \mathbb{R}^{N' \times N''} ,
 \forall c \in \mathcal{S}^{N''} ,\forall Q \in \mathbb{O}_{N'} ,\forall B \in \mathcal{S}^{N' - 1} }
 \right.\},
$$

and for given $ Q \in \mathbb{O}_{N'}$

\begin{align*}
&\mathcal{S}_{N' - 1} (Q) = \{ \left( {\begin{matrix}
   {Q\left( {\begin{matrix}
   0 & 0  \\
   0 & B  \\
 \end{matrix} } \right)Q^T } & {Qb}  \\
   {b^TQ^T} & c  \\
 \end{matrix} } \right) \in \mathcal{S}^N \left| {\forall b \in \mathbb{R}^{N' \times N''} ,
 \forall c \in \mathcal{S}^{N''} ,\forall B \in \mathcal{S}^{N' - 1} }
 \right.\}\\
 = & \{ \left( {\begin{matrix}
   Q & 0  \\
   0 & {I_{N''} }  \\
 \end{matrix} } \right)\left( {\begin{matrix}
   {\left( {\begin{matrix}
   0 & 0  \\
   0 & B  \\
 \end{matrix} } \right)} & b  \\
   {b^T } & c  \\
 \end{matrix} } \right)\left( {\begin{matrix}
   {Q^T } & 0  \\
   0 & {I_{N''} }  \\
\end{matrix} } \right)\left| {\forall b \in \mathbb{R}^{N' \times N''} ,
\forall c \in \mathcal{S}^{N''} ,\forall B \in \mathcal{S}^{N' - 1}
} \right.\}.
\end{align*}
Therefore $ \mathcal{S}_{N' - 1} (Q) \subset \mathcal{S}_{N' - 1}
\subset \mathcal{S}^N $. For any $(p',x'')$ fixed and $ Q \in
\mathbb{O}_{N'}$, $(A,p'',u,x') \in \mathcal{S}_{N' - 1} (Q)\times
\mathbb{R}^{N''} \times \mathbb{R} \times \mathbb{R}^{N'} $, we set
$$
X_F^*  = ((F^{ab}(A, p,u,x) ),F^{p_{N' + 1} } , \cdots ,F^{p_N }
,F^u ,F^{x_1 } , \cdots ,F^{x_{N'} } )
$$
as a vector in $\mathcal{S}^N  \times \mathbb{R}^{N''} \times
\mathbb{R} \times \mathbb{R}^{N'} $. Set
\begin{equation}
\Gamma _{X_F^* }^ \bot   = \{ \widetilde X \in \mathcal{S}_{N' - 1}
(Q) \times \mathbb{R}^{N''}  \times \mathbb{R} \times
\mathbb{R}^{N'} |\left\langle {\widetilde X, X_F^* } \right\rangle =
0\}.
\end{equation}

Let $ B \in \mathcal{S}^{N' - 1} $, $ A = B^{ - 1} $, and

\begin{center}
$\widetilde B = \left( {\begin{matrix}
   0 & 0  \\
   0 & B  \\
 \end{matrix} } \right)
$,
$ \widetilde A = \left( {\begin{matrix}
   0 & 0  \\
   0 & A  \\
 \end{matrix} } \right)
$.
\end{center}

For any given $ Q \in \mathbb{O}_{N'}$ and $ \widetilde X = ((X_{ab}
),(X_\alpha  ),Y,(Z_i )) \in \mathcal{S}_{N'-1}(Q) \times
\mathbb{R}^{N''} \times \mathbb{R} \times \mathbb{R}^{N'} $, we
define a quadratic form

\begin{align}
 Q^* (\widetilde X,\widetilde X)=&\sum\limits_{a,b,c,d = 1}^N {F^{ab,cd} X_{ab} X_{cd} }  +
2\sum\limits_{a,b = 1}^N {\sum\limits_{k,l = 1}^{N'} {F^{ab} [
{Q\widetilde A Q^T}]_{kl} X_{ka} X_{lb} } } \notag \\
&+ 2\sum\limits_{a,b = 1}^N {\sum\limits_{\alpha = N' + 1}^N
{F^{ab,p_\alpha  } X_{ab} X_\alpha  } } + 2\sum\limits_{a,b = 1}^N
{F^{ab,u} X_{ab} Y}  + 2\sum\limits_{a,b = 1}^N
{\sum\limits_{i = 1}^{N'} {F^{ab,x_i } X_{ab} Z_i } } \\
&+ \sum\limits_{\alpha ,\beta  = N' + 1}^N {F^{p_\alpha  ,p_\beta  }
X_\alpha  X_\beta  } + 2\sum\limits_{\alpha  = N' + 1}^N
{F^{p_\alpha  ,u} X_\alpha Y} + 2\sum\limits_{\alpha  = N' + 1}^N
{\sum\limits_{i = 1}^{N'}{F^{p_\alpha  ,x_i } X_\alpha  Z_i } } \notag \\
&+ F^{u,u} Y^2  + 2\sum\limits_{i = 1}^{N'} {F^{u,x_i } YZ_i }  +
\sum\limits_{i,j = 1}^{N'} {F^{x_i ,x_j } Z_i Z_j }, \notag
\end{align}

where the derivative functions of $F$ are evaluated at
$(\left(
{\begin{matrix}
   {Q\widetilde B Q^T } & {Qb}  \\
   {b^TQ^T} & c  \\
 \end{matrix} } \right),p,u,x)
$.

From lemma 3.3, we can get

\begin{lemma}
If $F$ satisfies condition (1.4), then for each $(p',x'')$
\begin{align}
F( {\left( {\begin{matrix}
   0 & b  \\
   b^T & c  \\
 \end{matrix} } \right)}, p, u, x) \text{ is locally convex in }
 (c,p'',u,x'), \text{ and }Q^* (\widetilde X,\widetilde X) \geqslant 0, \forall \widetilde X
\in \Gamma _{X_F^* }^ \bot,
\end{align}
where $Q^*$ is defined in (4.2).
\end{lemma}
Proof. Taking $\varepsilon >0$ small enough such that $a={Q\left(
{\begin{matrix}
   \varepsilon & 0  \\
   0 & B+\varepsilon I_{N'-1}  \\
 \end{matrix} } \right)Q^T}$ is invertible, and using (3.13), where $ \widetilde X = ((X_{ab}
),(X_\alpha  ),Y,(Z_i ))\in \Gamma _{X_F^* }^ \bot $, then we can
obtain (4.3) when $\varepsilon \to 0$.

Theorem 1.2 is a direct consequence of the following theorem and
Lemma 4.1.
\begin{theorem}
Suppose $\Omega$ is a domain in $\mathbb{R}^N  = \mathbb{R}^{N'}
\times \mathbb{R}^{N'' } $ and $ F(A,p,u,x) \in C^{2,1}
(\mathcal{S}^N \times \mathbb{R}^N \times \mathbb{R} \times \Omega
)$ satisfies (1.2) and  (1.4). Let $u \in C^{3,1} (\Omega) $ is a
partial convex solution of (1.1). If $ W(x)=(u_{ij}(x))_{N' \times
N'}$ attains minimum rank $l$ at certain point $z_0 \in \Omega$,
then there is a neighborhood $\mathcal {O}$ of $z_0$ and a positive
constant $C$ independent of $\phi$ (defined in (3.2)), such that
\begin{equation}
\sum\limits_{a,b = 1}^N {F^{ab} \phi _{ab} } \leqslant C(\phi  +
\left| {\nabla \phi } \right|), \quad \forall  x \in \mathcal {O}.
\end{equation}
In turn, $W(x)$is of constant rank in $\mathcal {O}$.
\end{theorem}

\textbf{Proof of Theorem 4.2}. Let $u \in C^{3,1}(\Omega)$ be a
partial convex solution of equation (1.1) and $ W(x)=(u_{ij}(x))_{N'
\times N'}$. For each $z_0 \in \Omega$ where $W$ attains minimal
rank $l$. We may assume $l \leqslant N'-1$, otherwise there is
nothing to prove. As in the previous section, we pick an open
neighborhood $\mathcal {O}$ of $z_0$, and for any $x \in \mathcal
{O}$, let $ G = \{ N' - l+1, \cdots ,N'\} $ and $ B = \{ 1,2, \cdots
,N'- l\} $ which means good terms and bad ones in indices for
eigenvalues of $W(x)$ respectively.

Setting $\phi$ as (3.2), then we see from Proposition 3.1 that
$$
\phi \in C^{1,1}(\mathcal {O}) ,\quad \phi(x) \geqslant 0, \quad
\phi(z_0) = 0,
$$
and there is a constant $C > 0$ such that for all $x \in \mathcal
{O}$,
\begin{equation}
\frac{1} {C}\sigma _1 (B)(x) \leqslant \phi (x) \leqslant C\sigma _1
(B)(x), \quad \frac{1} {C}\sigma _1 (B)(x) \leqslant \sigma _{l + 1}
(W(x)) \leqslant C\sigma _1 (B)(x).
\end{equation}
We shall fix a point $x \in \mathcal {O}$ and prove (4.4) at $x$.
For each $x \in \mathcal {O}$ fixed, we rotate coordinate $e_1,
\cdots, e_{N'}$ such that the matrix ${u_{ij}},i,j=1, \cdots, N'$ is
diagonal and without loss of generality we assume $ u_{11} \leqslant
u_{22} \leqslant \cdots \leqslant u_{N'N'} $. Then there is a
positive constant $C > 0$ depending only on $\left\| u
\right\|_{C^{3,1} }$ and $\mathcal {O}$, such that $ u_{N'N'}
\geqslant \cdots \geqslant u_{N' - l+1N' - l+1} \geqslant C > 0 $
for all $x \in \mathcal {O}$. Without confusion we will also simply
denote $ B = \{ u_{11} , \cdots ,u_{N'- lN'- l} \} $ and $ G = \{
u_{N' - l+1N' - l+1} , \cdots ,u_{N'N'} \} $. Note that for any
$\delta > 0$, we may choose $\mathcal {O}$ small enough such that
$u_{jj} < \delta$ for all $j \in B$ and $x \in \mathcal {O}$.

Again, as in section 3, we will avoid to deal with $\sigma_{l+1}(W)
= 0$. By considering $W_\varepsilon  = W + \varepsilon I$, and
$u_\varepsilon(x)=u(x)+ \frac{\varepsilon } {2}\left| x' \right|^2 $
for $\varepsilon >0$ sufficient small. Thus $u_\varepsilon(x)$
satisfies equation
\begin{equation}
F(D^2 u_\varepsilon,Du_\varepsilon,u_\varepsilon,x) =
R_\varepsilon(x),
\end{equation}
where $R_\varepsilon(x)=F(D^2
u_\varepsilon,Du_\varepsilon,u_\varepsilon,x)-F(D^2 u,Du,u,x)$.
Since $u \in C^{3,1}$; we have,
\begin{equation}
\left| {R_\varepsilon  (x)} \right| \leqslant C\varepsilon , \quad
\left| {\nabla R_\varepsilon  (x)} \right| \leqslant C\varepsilon ,
\quad \left| {\nabla ^2 R_\varepsilon  (x)} \right| \leqslant
C\varepsilon ,\quad \forall x \in \mathcal {O}.
\end{equation}

We will work on equation (4.6) to obtain differential inequality
(4.4) for $\phi_\varepsilon (x)$ defined in (3.3) with constant
$C_1$, $C_2$ independent of $\varepsilon$. Theorem 4.2 would follow
by letting $ \varepsilon \to 0$. In the following, we may as well
omit the subindex $\varepsilon$ for convenience.

We note that by (3.4), we have
$$
\varepsilon  \leqslant C\phi (x), \quad \forall x \in \mathcal {O},
$$
with $R(x)$ under control as follows,
\begin{equation}
\left|{D^jR_\varepsilon (x)} \right| \leqslant C\varepsilon, \text{
for all } j = 0, 1, 2, \text{ and for all } x \in \mathcal {O}.
\end{equation}
 Differentiate (4.6) one time in $x_i$ for $i \in B$, then we can get
$$
\sum\limits_{a,b = 1}^N {F^{ab} u_{abi} }  + \sum\limits_{a = 1}^N
{F^{p_a } u_{ai} }  + F^u u_i  + F^{x_i }  = O(\phi ),
$$
i.e.
\begin{align}
\sum\limits_{a,b = N' - l+1}^N {F^{ab} u_{abi} }  + \sum\limits_{a =
N' +1}^N {F^{p_a } u_{ai} }  + F^u u_i  + F^{x_i } =
O(\sum\limits_{i,j \in B} {\left| {\nabla u_{ij} } \right|}  + \phi
).
\end{align}
Differentiate (4.6) twice in $x_i$ for $i \in B$, then we obtain
\begin{align}
&\sum\limits_{a,b = 1}^N {F^{ab} u_{abii} }  + \sum\limits_{a,b =
1}^N {u_{abi} [\sum\limits_{c,d = 1}^N {F^{ab,cd} u_{cdi} }  +
\sum\limits_{c = 1}^N {F^{ab,p_c } u_{ci} }  + F^{ab,u} u_i  +
F^{ab,x_i } ]} \notag \\
& + \sum\limits_{a = 1}^N {F^{p_a } u_{aii} }  + \sum\limits_{a =
1}^N {u_{ai} [\sum\limits_{c,d = 1}^N {F^{p_a ,cd} u_{cdi} }  +
\sum\limits_{c = 1}^N {F^{p_a ,p_c } u_{ci} }  + F^{p_a ,u} u_i  + F^{p_a ,x_i } ]} \\
& + F^u u_{ii}  + u_i [\sum\limits_{c,d = 1}^N {F^{u,cd} u_{cdi} }
+ \sum\limits_{c = 1}^N {F^{u,p_c } u_{ci} }  + F^{u,u} u_i  + F^{u,x_i }  ] \notag \\
&+ \sum\limits_{c,d = 1}^N {F^{x_i ,cd} u_{cdi} }  + \sum\limits_{c
= 1}^N {F^{x_i ,p_c } u_{ci} }  + F^{x_i ,u} u_i  + F^{x_i ,x_i }  =
O(\phi ), \notag
\end{align}

i.e.
\begin{align}
&\sum\limits_{a,b = 1}^N {F^{ab} u_{abii} }  + \sum\limits_{a,b,c,d
= N' - l+1}^N {F^{ab,cd} u_{abi} u_{cdi} }  + 2\sum\limits_{a,b = N'
-l+1}^N {\sum\limits_{c = N' + 1}^N {F^{ab,p_c } u_{abi} u_{ci} } } \notag \\
& + 2\sum\limits_{a,b = N' - l+1}^N {F^{ab,u} u_{abi} u_i }  +
2\sum\limits_{a,b = N' - l+1}^N {F^{ab,x_i } u_{abi} }  +
\sum\limits_{a,c = N' + 1}^N {F^{p_a ,p_c } u_{ai} u_{ci} }\\
& + 2\sum\limits_{a = N' + 1}^N {F^{p_a ,u} u_{ai} u_i }  +
2\sum\limits_{a = N' + 1}^N {F^{p_a ,x_i } u_{ai} }  + F^{u,u} u_i
^2  + 2F^{u,x_i } u_i  + F^{x_i ,x_i } \notag \\
 = &O(\sum\limits_{i,j \in B} {\left| {\nabla u_{ij} } \right|}  + \phi
 ). \notag
\end{align}

So for each $i \in B$, let
\begin{align}
J_i= &\sum\limits_{a,b,c,d = N' - l+1}^N {F^{ab,cd} u_{abi} u_{cdi}
} + 2\sum\limits_{a,b = N' -l+1}^N {\sum\limits_{c = N' + 1}^N
{F^{ab,p_c } u_{abi} u_{ci} } } + 2\sum\limits_{a,b = N' - l+1}^N
{F^{ab,u} u_{abi} u_i }\notag \\
& + 2\sum\limits_{j \in G} {\frac{1} {{u_{jj} }}\sum\limits_{a,b =
N' - l+1}^N {F^{ab} u_{ija} u_{ijb} } }+2\sum\limits_{a,b = N' -
l+1}^N {F^{ab,x_i } u_{abi} }  +\sum\limits_{a,c = N' + 1}^N
{F^{p_a ,p_c } u_{ai} u_{ci} }\notag \\
& + 2\sum\limits_{a = N' + 1}^N {F^{p_a ,u} u_{ai} u_i }  +
2\sum\limits_{a = N' + 1}^N {F^{p_a ,x_i } u_{ai} }  + F^{u,u} u_i
^2  + 2F^{u,x_i } u_i  + F^{x_i ,x_i }.
\end{align}
Substitute (4.11) and (4.12) into (3.5), then we obtain
\begin{eqnarray}
&\sum\limits_{ab = 1}^N {F^{ab} \phi _{ab} }&  = -\sum\limits_{i \in
B} {[\sigma _l (G) + \frac{{\sigma _1 ^2 (B\left| i \right.) -
\sigma _2 (B\left| i \right.)}} {{\sigma _1 ^2 (B)}}]}J_i  \notag \\
&&- \frac{1}{{\sigma _1 ^3 (B)}}\sum\limits_{i \in B}
{\sum\limits_{ab = 1}^N {F^{ab} } [\sigma _1 (B)u_{iia}  - u_{ii}
\sum\limits_{j \in B} {u_{jja} } ]} [\sigma _1 (B)u_{iib}  - u_{ii}
\sum\limits_{j \in B}{u_{jjb} } ]\\
&&- \frac{1}{{\sigma _1 (B)}}\sum\limits_{\scriptstyle i,j \in B
\hfill \atop \scriptstyle i \ne j \hfill}  {\sum\limits_{ab = 1}^N
{F^{ab} u_{ija} u_{ijb} }} \notag \\
&&+ O(\sum\limits_{i,j \in B} {\left| {\nabla u_{ij} } \right|}  +
\phi ). \notag
\end{eqnarray}

By condition (1.4), since $u \in C^{3,1}$, so $F^{ab} \in C^{0,1}$.
For $ \overline {\mathcal {O}}  \subset \Omega $, there exists a
constant $ \delta_0 > 0$, such that
\begin{equation}
(F^{ab} ) \geqslant \delta _0 I_N, \quad \forall x \in \mathcal {O}.
\end{equation}

\textbf{Case(i): } $l=0$. Then $G= \emptyset$ and
\begin{align}
J_i= &\sum\limits_{a,b,c,d = N' + 1}^N {F^{ab,cd}(D^2u,Du,u,x)
u_{abi} u_{cdi} } + 2\sum\limits_{a,b = N' + 1}^N {\sum\limits_{c =
N' + 1}^N {F^{ab,p_c } u_{abi} u_{ci} } } \notag \\
& + 2\sum\limits_{a,b =N' + 1}^N {F^{ab,u} u_{abi} u_i }
+2\sum\limits_{a,b = N' + 1}^N {F^{ab,x_i } u_{abi} }
+\sum\limits_{a,c = N' + 1}^N {F^{p_a ,p_c } u_{ai} u_{ci} } \\
&+ 2\sum\limits_{a = N' + 1}^N {F^{p_a ,u} u_{ai} u_i }
+2\sum\limits_{a = N' + 1}^N {F^{p_a ,x_i } u_{ai} }  + F^{u,u} u_i
^2  + 2F^{u,x_i } u_i  + F^{x_i ,x_i }, \notag
\end{align}
where all the derivative functions of $F$ are evaluated at
$(D^2u,Du,u,x)$. Since $F \in C^{2,1}$ and $ \left\| {\left. {W(x)}
\right\|_{C^0 } } \right. = O(\phi ) $, by Taylor formula and
condition (4.3), we can get
\begin{align}
J_i = & O(\phi)+\sum\limits_{a,b,c,d = N' + 1}^N {F^{ab,cd} u_{abi}
u_{cdi} } + 2\sum\limits_{a,b = N' + 1}^N {\sum\limits_{c =
N' + 1}^N {F^{ab,p_c } u_{abi} u_{ci} } } \notag \\
& + 2\sum\limits_{a,b =N' + 1}^N {F^{ab,u} u_{abi} u_i }
+2\sum\limits_{a,b = N' + 1}^N {F^{ab,x_i } u_{abi} }
+\sum\limits_{a,c = N' + 1}^N {F^{p_a ,p_c } u_{ai} u_{ci} } \notag \\
&+ 2\sum\limits_{a = N' + 1}^N {F^{p_a ,u} u_{ai} u_i }
+2\sum\limits_{a = N' + 1}^N {F^{p_a ,x_i } u_{ai} }  + F^{u,u} u_i
^2  + 2F^{u,x_i } u_i  + F^{x_i ,x_i } \notag \\
\geqslant& -C\phi,
\end{align}
where all the derivative functions of $F$ are evaluated at
$(\left(
{\begin{matrix}
   0 & {(u_{k \alpha})}  \\
   {(u_{\alpha k})} & {(u_{\alpha \beta})}  \\
 \end{matrix} } \right),p,u,x)
$.

\textbf{Case(ii): } $1 \leqslant l \leqslant N'-1$ \\
 Now we set $ X_{ab}  = 0 $ for $a \in B $ or $b \in B $,
\begin{equation}
X_{N'N'}  =u_{N'N'i}  - \frac{1} {{F^{N'N'} }}[\sum\limits_{a,b = N'
- l+1}^N {F^{ab} u_{abi} }  + \sum\limits_{a = N' +1}^N {F^{p_a }
u_{ai} } + F^u u_i  + F^{x_i }],
\end{equation}
 $ X_{ab}  = u_{abi} $ otherwise, $Y=u_i$ and $Z_k=\delta_{ki}$.
 We can verify that $(X_{ab}) \in S_{N'-1}(I_{N'})$ and $ \widetilde X = ((X_{ab}
),(X_\alpha  ),Y,(Z_i )) \in \Gamma _{X_F^* }^ \bot $. Again by
condition (4.3), we infer that
\begin{align}
J_i \geqslant -C(\sum\limits_{i,j \in B} {\left| {\nabla u_{ij} }
\right|}  + \phi ),
\end{align}
since $C>\sigma _l (G) + \frac{{\sigma _1 ^2 (B\left| i \right.) -
\sigma _2 (B\left| i \right.)}} {{\sigma _1 ^2 (B)}}>0$, thus we
obtain
\begin{eqnarray}
\sum\limits_{a,b = 1}^N {F^{ab} \phi _{ab} } &\leqslant&
C(\sum\limits_{i,j \in B}{\left| {\nabla u_{ij} } \right|}  + \phi ) \notag \\
&&- \frac{1}{{\sigma _1 ^3 (B)}}\sum\limits_{i \in B}
{\sum\limits_{a,b = 1}^N {F^{ab} } [\sigma _1 (B)u_{iia}  - u_{ii}
\sum\limits_{j \in B} {u_{jja} } ]} [\sigma _1 (B)u_{iib}  - u_{ii}
\sum\limits_{j \in B}{u_{jjb} } ] \notag \\
&&- \frac{1}{{\sigma _1 (B)}}\sum\limits_{\scriptstyle i,j \in B
\hfill \atop \scriptstyle i \ne j \hfill}  {\sum\limits_{a,b = 1}^N
{F^{ab} u_{ija} u_{ijb} } } \notag \\
&\leqslant& C(\sum\limits_{i,j \in B} {\left| {\nabla u_{ij} }
\right|}  + \phi ) - \frac{\delta_0}{{\sigma _1 ^3
(B)}}\sum\limits_{i \in B} {\sum\limits_{a = 1}^N \widetilde V_{i
a}^2}- \frac{\delta_0}{{\sigma _1 (B)}}\sum\limits_{\scriptstyle i,j
\in B \hfill \atop \scriptstyle i \ne j \hfill}  {\sum\limits_{a =
1}^N  u_{ija}^2 },
\end{eqnarray}
where $\widetilde V_{i a}=\sigma _1 (B)u_{iia}  - u_{ii}
\sum\limits_{j \in B} {u_{jja} } $. Referring to Lemma 3.3 in [2],
we can control the term $\sum\limits_{i,j \in B} {\left| {\nabla
u_{ij} } \right|}$ by the rest terms on the right hand side in
(4.19) and $ \phi+\left| {\nabla \phi } \right| $ where
\begin{equation}
\phi_a=O(\phi)+ \sum\limits_{i \in B}[\sigma _l (G) + \frac{{\sigma
_1 ^2 (B\left| i \right.) - \sigma _2 (B\left| i \right.)}} {{\sigma
_1 ^2 (B)}}]u_{iia}.
\end{equation}
So there exist positive constants $C_1$,$C_2$ independent of
$\varepsilon$, such that
\begin{equation}
\sum\limits_{ab = 1}^N {F^{ab} \phi _{ab} } \leqslant C_1 (\phi  +
\left| {\nabla \phi } \right|)-C_2\sum\limits_{i,j \in B} {\left|
{\nabla u_{ij} } \right|}, \quad \forall  x \in \mathcal {O}.
\end{equation}
Taking $\varepsilon \to 0$, (4.19) is proved for $u$. By the Strong
Maximum Principle, $\phi (x) \equiv 0 $ in $\mathcal {O}$; and $W$
is of constant rank in $\mathcal {O}$. The proof of Theorem 4.2 is
completed.
\begin{remark}
In the above proof, we have used a weak condition (4.3). Also we can
directly use the condition (1.4), i.e.(3.13). We set $ X_{ab}  = 0 $
for $a \in B $ or $b \in B $, $ X_{ab}  = u_{abi} $ otherwise,
$Y=u_i$ and $Z_k=\delta_{ki}$. Then we have $ \widetilde X =
((X_{ab} ),(X_\alpha ),Y,(Z_i )) \in \mathcal{S}^N \times
\mathbb{R}^{N''} \times \mathbb{R} \times \mathbb{R}^{N'}$, and by
(3.13), $J_i \geqslant 0$ for every $i \in B$. So (4.19) holds. As
above, Theorem 4.2 holds.
\end{remark}
\begin{remark}
In particular, for $N'=1$, we only need the following structure
condition
\begin{equation}
 F(\left( {\begin{matrix}
   0 & b  \\
   {b^T } & c  \\
\end{matrix} } \right),p',p'',u,x',x'') \text{ is locally convex in }
(c,p'',u,x'),
\end{equation}
then we have $(u_{ij})_{N' \times N'}$ is of constant rank in
$\Omega$. Since when $N'=1$, so the minimum rank $l$ has only two
cases: $l=1$ and $l=0$. If $l=1$ we are done; and if $l=0$, (4.16)
and (4.19) holds by condition (4.22). Then the result holds as the
proof of Theorem 4.2.
\end{remark}

\section{the proof of theorem 1.5}

In this section we give the proof of Theorem 1.5. It is similar to
the proof of Theorem 1.2 only some minor modifications.

Following the notations of Theorem 1.5, suppose
$W(x,t_0)=(u_{ij}(x,t_0))_{N' \times N'}$ attains minimal rank
$l=l(t_0)$ at some point $z_0 \in \Omega$. We may assume $l\leqslant
N'-1$, otherwise there is nothing to prove. As in the section 4,
there is a neighborhood $\mathcal {O}\times (t_0-\delta,
t_0+\delta]$ of $(z_0, t_0)$ instead of $\mathcal {O}$, such that  $
u_{N'N'} \geqslant \cdots \geqslant u_{N' - l+1N' - l+1} \geqslant C
> 0 $ for all $(x,t) \in \mathcal {O}\times (t_0-\delta,
t_0+\delta]$, and we can denote $ B = \{ u_{11} , \cdots ,u_{N'-
lN'- l} \} $ and $ G = \{ u_{N' - l+1N' - l+1} , \cdots ,u_{N'N'} \}
$. If $t_0=T$, the neighborhood should be $\mathcal {O}\times
(t_0-\delta, t_0]$.

Setting $\phi$ as (3.2) (where $W(x,t)$ instead of $W(x)$), then we
see from Proposition 3.1 that
$$
\phi \in C^{1,1}(\mathcal {O}\times (t_0-\delta, t_0+\delta]) ,\quad
\phi(x,t) \geqslant 0, \quad \phi(z_0,t_0) = 0,
$$
Also when we choose $\mathcal {O}$ and $\delta>0$ small enough, the
corresponding (3.4), (4.5) and (4.8) hold. Then Theorem 1.5 is a
consequence of the following theorem and the method of continuity.

\begin{theorem}
Suppose $\Omega$ is a domain in $\mathbb{R}^N  = \mathbb{R}^{N'}
\times \mathbb{R}^{N'' } $ and $ F(A,p,u,x,t) \in C^{2,1}
(\mathcal{S}^N \times \mathbb{R}^N \times \mathbb{R} \times \Omega
\times (0,T])$ satisfies (1.2) for each $t$ and  (1.8). Let $u \in
C^{3,1}$ is a partial convex solution of (1.9). For each $t_0 \in
(0,T]$, if $ W(x,t_0)=(u_{ij}(x,t_0))_{N' \times N'}$ attains
minimum rank $l$ at some point $z_0 \in \Omega$, then there is a
neighborhood $\mathcal {O}\times (t_0-\delta, t_0+\delta] $ of
$(z_0, t_0)$ as above and a positive constant $C$ independent of
$\phi$ (defined in (3.2)), such that
\begin{equation}
\sum\limits_{ab = 1}^N {F^{ab} \phi _{ab}(x,t)-\phi_t(x,t) }
\leqslant C(\phi (x,t) + \left| {\nabla \phi(x,t) } \right|), \quad
\forall  (x,t) \in \mathcal {O} \times (t_0-\delta, t_0+\delta].
\end{equation}
In turn, $W(x,t)$ has constant rank $l$ in $\mathcal {O} \times
(t_0-\delta, t_0]$, where $l=l(t_0)$.
\end{theorem}

\textbf{Proof of Theorem 5.1}.The proof is similar to the proof of
Theorem 4.2,  so we only give the main process of the proof.

With $u_t = F(D^2u,Du, u, x, t)$, using the same notations as above
and the proof of Theorem 4.2, we have $u_\varepsilon(x,t)=u(x,t)+
\frac{\varepsilon } {2}\left| x' \right|^2 $ for $\varepsilon >0$
sufficient small. Thus $u_\varepsilon(x,t)$ satisfies equation
\begin{equation}
(u_\varepsilon)_t
=F(D^2u_\varepsilon,Du_\varepsilon,u_\varepsilon,x,t) -
R_\varepsilon(x,t),
\end{equation}
where $R_\varepsilon(x,t)=F(D^2
u_\varepsilon,Du_\varepsilon,u_\varepsilon,x,t)-F(D^2 u,Du,u,x,t)$.

As in the proof of Theorem 4.2, we omit the subindex $\varepsilon$.

Differentiate (5.2) one time in $x_i$ for $i \in B$, then we can get
$$
\sum\limits_{a,b = 1}^N {F^{ab} u_{abi} }  + \sum\limits_{a = 1}^N
{F^{p_a } u_{ai} }  + F^u u_i  + F^{x_i }  = O(\phi )+u_{i,t} \text{
,}
$$
i.e.
\begin{align}
\sum\limits_{a,b = N' - l+1}^N {F^{ab} u_{abi} }  + \sum\limits_{a =
N' +1}^N {F^{p_a } u_{ai} }  + F^u u_i  + F^{x_i } =
O(\sum\limits_{i,j \in B} {\left| {\nabla u_{ij} } \right|}  + \phi
)+u_{i,t} \text{ .}
\end{align}
Differentiate (5.2) twice in $x_i$ for $i \in B$, then we can get
\begin{align}
&\sum\limits_{a,b = 1}^N {F^{ab} u_{abii} }  + \sum\limits_{a,b =
1}^N {u_{abi} [\sum\limits_{c,d = 1}^N {F^{ab,cd} u_{cdi} }  +
\sum\limits_{c = 1}^N {F^{ab,p_c } u_{ci} }  + F^{ab,u} u_i  +
F^{ab,x_i } ]} \notag \\
& + \sum\limits_{a = 1}^N {F^{p_a } u_{aii} }  + \sum\limits_{a =
1}^N {u_{ai} [\sum\limits_{c,d = 1}^N {F^{p_a ,cd} u_{cdi} }  +
\sum\limits_{c = 1}^N {F^{p_a ,p_c } u_{ci} }  + F^{p_a ,u} u_i  + F^{p_a ,x_i } ]} \notag \\
& + F^u u_{ii}  + u_i [\sum\limits_{c,d = 1}^N {F^{u,cd} u_{cdi} }
+ \sum\limits_{c = 1}^N {F^{u,p_c } u_{ci} }  + F^{u,u} u_i  + F^{u,x_i }  ] \\
&+ \sum\limits_{c,d = 1}^N {F^{x_i ,cd} u_{cdi} }  + \sum\limits_{c
= 1}^N {F^{x_i ,p_c } u_{ci} }  + F^{x_i ,u} u_i  + F^{x_i ,x_i }  =
O(\phi )+u_{ii,t} \text{ ,} \notag
\end{align}
i.e.
\begin{align}
&\sum\limits_{a,b = 1}^N {F^{ab} u_{abii} }  + \sum\limits_{a,b,c,d
= N' - l+1}^N {F^{ab,cd} u_{abi} u_{cdi} }  + 2\sum\limits_{a,b = N'
-l+1}^N {\sum\limits_{c = N' + 1}^N {F^{ab,p_c } u_{abi} u_{ci} } } \notag \\
& + 2\sum\limits_{a,b = N' - l+1}^N {F^{ab,u} u_{abi} u_i }  +
2\sum\limits_{a,b = N' - l+1}^N {F^{ab,x_i } u_{abi} }  +
\sum\limits_{a,c = N' + 1}^N {F^{p_a ,p_c } u_{ai} u_{ci} }\\
& + 2\sum\limits_{a = N' + 1}^N {F^{p_a ,u} u_{ai} u_i }  +
2\sum\limits_{a = N' + 1}^N {F^{p_a ,x_i } u_{ai} }  + F^{u,u} u_i
^2  + 2F^{u,x_i } u_i  + F^{x_i ,x_i } \notag \\
= &O(\sum\limits_{i,j \in B} {\left| {\nabla u_{ij} } \right|}  +
\phi)+u_{ii,t} \text{ .}\notag
\end{align}
We denote that
$$ \phi _t  = \sum\limits_{i,j = 1}^{N'}
{\frac{{\partial \phi }} {{\partial u_{ij} }}u_{ij,t}  =
\sum\limits_{i = 1}^{N'} {\frac{{\partial \phi }} {{\partial u_{ii}
}}u_{ii,t} } } \text{ ,}
$$
so we can obtain from (3.5), (4.12) and (5.5),
\begin{eqnarray}
&&\sum\limits_{ab = 1}^N {F^{ab} \phi _{ab}(x,t)-\phi_t(x,t)} =
-\sum\limits_{i \in B} {[\sigma _l (G) + \frac{{\sigma _1 ^2
(B\left| i \right.) - \sigma _2 (B\left| i \right.)}} {{\sigma _1 ^2 (B)}}]}J_i  \notag \\
&&- \frac{1}{{\sigma _1 ^3 (B)}}\sum\limits_{i \in B}
{\sum\limits_{ab = 1}^N {F^{ab} } [\sigma _1 (B)u_{iia}  - u_{ii}
\sum\limits_{j \in B} {u_{jja} } ]} [\sigma _1 (B)u_{iib}  - u_{ii}
\sum\limits_{j \in B}{u_{jjb} } ]\\
&&- \frac{1}{{\sigma _1 (B)}}\sum\limits_{\scriptstyle i,j \in B
\hfill \atop \scriptstyle i \ne j \hfill}  {\sum\limits_{ab = 1}^N
{F^{ab} u_{ija} u_{ijb} }} \notag \\
&&+ O(\sum\limits_{i,j \in B} {\left| {\nabla u_{ij} } \right|}  +
\phi ) .\notag
\end{eqnarray}
Now the right hand side of (5.6) is the same as the right hand side
of (4.13). From Remark 4.3, we set $ X_{ab}  = 0 $ for $a \in B $ or
$b \in B $, $ X_{ab}  = u_{abi} $ otherwise, $Y=u_i$ and
$Z_k=\delta_{ki}$. Then we have $ \widetilde X = ((X_{ab}),(X_\alpha
),Y,(Z_i )) \in \mathcal{S}^N \times \mathbb{R}^{N''}  \times
\mathbb{R} \times \mathbb{R}^{N'}$, and by (3.13), $J_i \geqslant 0$
for every $i \in B$. So (4.19) holds. A similar analysis as in the
proof of Theorem 4.2 for the right hand side of equation (4.19)
yields
\begin{equation}
\sum\limits_{ab = 1}^N {F^{ab} \phi _{ab}(x,t)-\phi_t(x,t) }
\leqslant C_1(\phi (x,t) + \left| {\nabla \phi(x,t) }
\right|)-C_2\sum\limits_{i,j \in B} {\left| {\nabla u_{ij} }
\right|},
\end{equation}
where the positive constants $C_1$,$C_2$ independent of
$\varepsilon$, and $(x,t) \in \mathcal {O} \times (t_0-\delta,
t_0+\delta]$. Then $W(x,t)$ has a constant rank $l$ for each $(x,t)
\in \mathcal {O} \times (t_0-\delta, t_0]$ by the Strong Maximum
Principle for parabolic equations. Theorem 5.1 holds.

\section{discussion of structure condition }

In this section, we discuss the condition (4.3) and (1.4).

For any given $ Q \in \mathbb{O}_{N'}$, we define
$$
\widetilde{F_Q }(A,b,c,p'',u,x') = F(\left( {\begin{matrix}
   {Q\left( {\begin{matrix}
   0 & 0  \\
   0 & {A^{ - 1} }  \\
 \end{matrix} } \right)Q^T } & {Q\left( {\begin{matrix}
   0 & 0  \\
   0 & {A^{ - 1} }  \\
 \end{matrix} } \right)b}  \\
   {b^T \left( {\begin{matrix}
   0 & 0  \\
   0 & {A^{ - 1} }  \\
 \end{matrix} } \right)Q^T } & {c + b^T \left( {\begin{matrix}
   0 & 0  \\
   0 & {A^{ - 1} }  \\
 \end{matrix} } \right)b}  \\
 \end{matrix} } \right),p,u,x),
$$
for $ (A,b,c,p'',u,x') \in \mathcal{S}_+^{N' - 1}  \times
\mathbb{R}^{N' \times N''} \times \mathcal{S}^{N''}  \times
\mathbb{R}^{N''} \times \mathbb{R} \times \mathbb{R}^{N'} $ and
fixed $ (p',x'') \in \mathbb{R}^{N'}  \times \mathbb{R}^{N''} $.
Condition (1.4) implies the following condition
\begin{equation}
\widetilde {F}_Q (A,b,c,p'',u,x') \text{ is locally convex in }
(A,b,c,p'',u,x'),
\end{equation}
for any fixed $N' \times N'$ orthogonal matrix $Q$.

\begin{proposition}
Let $Q \in \mathbb{O}_{N'}$. The condition (6.1)is equivalent to
\begin{equation}
Q^* (\widetilde X,\widetilde X) \geqslant 0,
\end{equation}
for any $ \widetilde X = ((X_{ab} ),(X_\alpha  ),Y,(Z_i )) \in
\mathcal{S}_{N'-1}(Q) \times \mathbb{R}^{N''} \times \mathbb{R}
\times \mathbb{R}^{N'} $, where $Q^*$ is defined in (4.2).
\end{proposition}
\textbf{Proof}. By approximating, Proposition 6.1 holds from Lemma
3.2.

\begin{remark}
Condition (1.4) is equivalent to (3.13), and (1.4) implies (6.1) for
any fixed $N' \times N'$ orthogonal matrix $Q$. Condition (6.1) is
equivalent to (6.2), and Lemma 4.1 is a consequence of Proposition
6.1. And condition (6.1) is weaker than condition (1.4).
\end{remark}

There is a class of functions which satisfy (1.4). Through a direct
calculation and using (3.13), we can get
\begin{proposition}
If $g$ is a non-decreasing and convex function and $F_1$, $\cdots$,
$F_m$ satisfy condition (1.4), then $F = g(F_1,\cdots, F_m)$ also
satisfies condition (1.4). In particular,  if $F_1$ and $F_2$ are in
the class, so are $F_1 +F_2$ and $F_1^\alpha$( where $F_1 > 0$) for
any $\alpha \geqslant 1$.
\end{proposition}

\begin{remark}
This paper was finished in April 2009, and B. Bian and P. Guan give
a better structural condition ( an equivalent condition of (4.3)) in
their paper "A Structural Condition for Microscopic Convexity
Principle", which appears in Discrete and Continuous Dynamical
Systems, Volume 28, Number 2, 2010, pp. 789-807.
\end{remark}

\end{document}